\documentclass[aop,MSNbibl,seceqn,citesort,dvips]{arximspdf}

\doi{10.1214/11-AOP687}
\volume{40}
\issue{6}
\pubyear{2012}
\firstpage{2539}
\lastpage{2588}

\makeatletter

\newtheorem{theorem}{Theorem}[section]
\newtheorem{pro}[theorem]{Proposition}
\newtheorem{lem}[theorem]{Lemma}
\newtheorem{cor}[theorem]{Corollary}

\newproclaim{defin}[theorem]{Definition}
\newproclaim{exa}{Example}[section]
\newproclaim{rem}[theorem]{Remark}

\newcommand{\ra}{\rightarrow}

\makeatother

\begin{document}
\begin{frontmatter}

\title{Random walks driven by low moment measures}
\runtitle{Low moment random walks}

\begin{aug}
\author[A]{\fnms{Alexander} \snm{Bendikov}\thanksref{t2}\ead[label=e1]{Alexander.Bendikov@math.uni.wroc.pl}} and
\author[B]{\fnms{Laurent} \snm{Saloff-Coste}\corref{}\thanksref{t3}\ead[label=e2]{lsc@math.cornell.edu}}
\runauthor{A. Bendikov and L. Saloff-Coste}
\affiliation{Wroc\l aw University and Cornell University}
\address[A]{Institute of Mathematics\\
Wroc\l aw University\\
Pl. Grundwaldzki 2/4\\
50-384 Wroc\l aw\\
Poland\\
\printead{e1}}
\address[B]{Department of Mathematics\\
Malott Hall\\
Cornell University\\
Ithaca, New York 14850-4201\\
USA\\
\printead{e2}}
\end{aug}

\thankstext{t2}{Supported in part by
Polish Government Scientific Research Fund, Grant NN201371736.}

\thankstext{t3}{Supported in part by NSF
Grant DMS-06-03886 and 104771.}

\received{\smonth{8} \syear{2010}}
\revised{\smonth{5} \syear{2011}}

%
\begin{abstract}
We study the decay of convolution powers of probability measures
without second moment but satisfying some weaker finite moment
condition. For any locally compact unimodular group $G$ and any
positive function $\varrho\dvtx G\ra[0,+\infty]$, we introduce a function
$\Phi_{G,\varrho}$ which describes the fastest possible decay of
$n\mapsto\phi^{(2n)}(e)$ when $\phi$ is a symmetric continuous
probability density such that $\int\varrho\phi$ is finite. We
estimate $\Phi_{G,\varrho}$ for a variety of groups $G$ and functions
$\varrho$. When $\varrho$ is of the form $\varrho=\rho\circ\delta$
with $\rho\dvtx[0,+\infty)\ra[0,+\infty)$, a fixed increasing function,
and $\delta\dvtx G\ra[0,+\infty)$, a natural word length measuring the
distance to the identity element in $G$, $\Phi_{G,\varrho}$ can be
thought of as a group invariant.
\end{abstract}

%

\end{frontmatter}

\section{Introduction} \label{sec-Intr}

Throughout this work, $G$ is a locally compact
unimodular group equipped with its Haar measure $\lambda$, and $L^p(G)=
L^p(G,\lambda)$, $1\le p\le\infty$,
is the space of (classes of) $p$ integrable measurable functions.
When convenient, we write $\lambda(dx)=dx$.

Sometimes, but not always, we will assume that $G$ is also compactly generated.
When that is the case, we let
$U$ be an open relatively compact set which is symmetric and contains a
compact generating neighborhood of the identity element $e$.
For any element $x$ in $G$, we set $|x|=\inf\{n\dvtx x\in U^n\}$ (with
the convention that $U^0=\{e\}$) and $V(n)=\lambda(U^n)$.
The function $V$ is called the volume growth function of the group $G$.
The rough behavior of both $x\mapsto|x|$ and $n\mapsto V(n)$
is essentially independent of the choice of $U$; for example, see
\cite{VSCC}.
The case when $G$ is
a finitely generated group equipped with a finite symmetric generating set,
and its counting measure is of course
included here, and the results we obtain are particularly
interesting in this case.

Given a Borel probability measure $\mu$ on $G$, we let $\mu^{(n)}$ be the
$n$-fold convolution power of $\mu$ and let $\check{\mu}$ be the
measure defined
by $\check{\mu}(A)=\mu(A^{-1})$ for any Borel set~$A$. Recall that $\mu^{(n)}$
is the law of the random walk driven by $\mu$ and started at $e$.
We call a measure symmetric if $\mu=\check{\mu}$.
Since $G$ is unimodular, we have $\check{\lambda}=\lambda$.
It follows that a measure having a density $\phi$ w.r.t.\
the Haar measure $\lambda$ is symmetric if and only if $\phi$ is
symmetric, that is, $\phi=\check{\phi}$ where $\check{\phi}(x)=\phi
(x^{-1})$; see, for example,~\cite{Dix2}, Exercise 5, page 89. Throughout
the paper, we denote by
$R_\phi$ the operator of convolution by the function $\phi\in L^2(G)$
on the right,\vspace*{1pt} that is, $R_\phi f=f*\phi$
(say, for compactly supported continuous function $f$).
When $\phi$ is in $L^1(G)$, $R_\phi$ also denotes the extension of this operator
to $L^2(G)$ [and, more generally, $L^p(G)$]. When
$\phi=\check{\phi}\in L^1(G)$, $R_\phi$ is a bounded
self-adjoint operator on $L^2(G)$.

\subsection{The decay of convolution powers}

A probability measure $\mu$ on a compactly generated group $G$ is said to
have finite second moment if $\mu(|\cdot|^2)<\infty$.
A fundamental result concerning symmetric random walks
on groups asserts that there exists a nonincreasing positive function
$\Phi_G$ such that, for any symmetric probability measure $d\mu=\phi
\,d\lambda$
with finite second moment and continuous density $\phi$
whose support contains a generating compact neighborhood of the identity,
we have
%
\begin{equation}\label{Phi}
\mu^{(2n)}(U)\simeq\phi^{(2n)}(e)\simeq\Phi_G(n);
\end{equation}
see~\cite{Gre,PSCstab}.
Here, $f(n)\simeq g(n)$ means that there are constants $c_i\in(0,\infty)$
such that, for all $n$, $c_1f(c_2n)\le g(n)\le c_3f(c_4 n)$.
Clearly, in the above estimates, the implied constants $c_i$ are
allowed to
depend on $\mu$ and $G$.

The following list provides some examples of explicit computation of
$\Phi_G$,
assuming that $G$ is compactly generated. More accurately, it is the
equivalence class of $\Phi_G$ under the equivalence relation $\simeq$
which is
computed.
\begin{itemize}
\item If $G$ is such that $V(n) \simeq n^D$, then $\Phi_G(n)\simeq n^{-D/2}$.
Every nilpotent group has these properties for some integer $D$; see
\cite{Varconv,VSCC} and the references therein.
\item If $G$ is polycyclic (or linear solvable) and
has exponential volume growth, then
$\Phi_G(n)\simeq\exp(-n^{1/3})$; see~\cite{Alex,PSCrank,Var-Sendai,VSCC}.
\item The group $G$ is nonamenable if and only if $\Phi_G(n)\simeq\exp(-n)$
(this is a formulation of Kesten's celebrated theorem regarding amenability
and random walks).

\item Let $M,N$ be two finitely generated groups, and let
$G$ be the wreath product $G= M\wr N=(\sum_{n\in N}M_n)\rtimes N$.
This is the semidirect product of $N$ with
of the direct sum of countably many
copies of $M$ indexed by $N$ where the action of $N$ is by index
translation; see, for example,~\cite{PSCwp} for a precise definition.
\begin{itemize}
\item Assume $N$ satisfies $V_N(n)\simeq n^d$ for some $d\ge1$ and $M$
is nontrivial. Then we have
\[
\Phi_G(n)\simeq\cases{
\exp\bigl(- n^{d/(d+2)}\bigr), &\quad if $M$ is finite;\vspace*{2pt}\cr
\exp\bigl(- [n^d(\log n)^2]^{1/(d+2)}\bigr), &\quad if $V_M(n)\simeq n^b, b\ge
1$;\vspace*{2pt}\cr
\exp\bigl(-n^{(d+1)/(d+3)}\bigr), &\quad if $M \in\mathcal P\mathcal E$,}
\]
where $\mathcal P\mathcal E$ stands for polycyclic
with exponential volume growth.
\item Assume that $N\in\mathcal P \mathcal E$ and $M$ is nontrivial,
finite or polycyclic. Then we have
\[
\Phi_G(n)\simeq\exp(-n(\log n)^{-2});
\]
\end{itemize}
see
\cite{Erschleriso,ErschlerPTRF,PSCwp,SCsurveydiff} for details and further
results.
\item Let $N=\mathbb Z^d$, $M$ be nontrivial, and $k\ge2$
be an integer. Set
$G= M\wr(M\wr(\cdots(M\wr N)\cdots))$ where $k$ successive wreath products
are taken. Then
\[
\Phi_G(n)\simeq\cases{
\exp(- n (\log_{k-1} n)^{-2/d}), &\quad if $M$ is finite,\vspace*{2pt}\cr
\exp\bigl(- n[(\log_{k-1} n)/\log_k n ]^{-2/d}\bigr), &\quad if $V_M(n)\simeq
n^b, b\ge1$.}
\]
Here, $\log_1(x)=\log(e+x)$ and $\log_k (x)= \log(e+\log_{k-1}(x))$,
$k\ge2$;
see~\cite{Erschleriso,ErschlerPTRF,SCsurveydiff}.
\end{itemize}
The article~\cite{SCnotices} gives an overview.
Many further behaviors are possible for the function $\Phi_G$, but a
complete classification of the possible behaviors is not known. In fact,
the very existence of such a classification seems highly unlikely, and
there are (uncountably) many amenable finitely generated groups $G$
for which the behavior of $\Phi_G$ is unknown.
Still, Definition~\ref{Phi} means that on any such group,
we know that all random walks driven by a symmetric measure with
generating support and
finite second moment have comparable probability of return.

This work focuses on the probability of return of random walks driven
by measures that may fail to have a finite second moment
but satisfy some finite moment condition. Namely,
consider a nonnegative, nondecreasing function
$\rho\dvtx[0,+\infty)\ra[0,+\infty)$. For any finitely generated
group $G$ equipped with a word length \mbox{$|\cdot|$} as above,
let $\rho_G$ be the function
\[
\rho_G\dvtx G\ra[0,+\infty),\qquad x\mapsto\rho_G(x)=\rho(|x|).
\]
We will abuse notation and write $\rho$ for $\rho_G$ when convenient.
We say that a probability measure $\mu$ on $G$ has finite $\rho_G$-moment
if
\[
\mu(\rho_G)=\sum_{g\in G}\rho_G(g)\mu(g)<\infty.
\]
Since we are mostly interested in measures without second moment, the following
are some of the natural choices for $\rho$:
\begin{itemize}
\item Small powers: $\rho_\alpha(t)= (1+t)^\alpha$, $\alpha\in(0,2)$.
\item Regularly varying functions of index $\alpha\in(0,2)$, for example,
\[
\rho(t)=(1+t)^\alpha[\log(e+t)]^\beta,\qquad \beta\in\mathbb R.
\]
\item Slowly varying increasing functions including:
\begin{itemize}
\item$\rho^{\exp}_{c,\alpha}(t)=\exp( c [\log(1+t)]^\alpha)$, $\alpha
\in(0,1)$ and $c>0$;
\item$\rho^{\log}_\alpha(t)= [\log(e+t)]^\alpha$, $\alpha\in(0,\infty)$.
\end{itemize}
\end{itemize}

We consider the following natural question.
What can be said about the decay of $\phi^{(2n)}(e)$ when
$d\mu=\phi \,d\lambda$ is a
symmetric measure having a finite $\rho_G$-moment for one of
the functions $\rho$ mentioned above?

\subsection{Group invariants associated with random walks and moment conditions}

In general, requiring that a symmetric measure $\mu$ has
a finite moment of some sort is not enough to determine
the behavior of the convolution powers
of that measure. The following definition introduces the notion
of ``fastest decay'' allowed by a given moment condition.
\begin{defin}[(Fastest decay under $\varrho$-moment)] \label{def1}
Let $G$ be a locally compact unimodular group.
Fix a measurable function $\varrho\dvtx G\ra[0,+\infty]$.
Fix a compact symmetric neighborhood\vspace*{1pt} $\Omega$ of $e$ in $G$
such that $\lambda(\Omega)\ge1$ and $\sup_{\Omega^2}\{\varrho\}>0$.
For $K> 1$, let $\mathcal S^{\Omega,K}_{G,\varrho}$ be the set of all
symmetric continuous
probability densities $\phi$ on $G$ with the properties that $\|\phi\|
_\infty\le K$ and
$\int\phi\varrho \,d\lambda\le K \sup_{\Omega^2}\{\varrho\}$. Set
\[
\Phi^{\Omega,K}_{G,\varrho}\dvtx n\mapsto\Phi^{\Omega,K}_{G,\varrho}(n):=
\inf\bigl\{
\phi^{(2n)}(e) \dvtx \phi\in\mathcal S^{\Omega,K}_{G,\varrho}\bigr\}.
\]
\end{defin}

In words, $\Phi^{\Omega,K}_{G,\varrho}$ provides the best lower bound
valid for all convolution powers of probability measures with density in
$\mathcal S ^{\Omega,K}_{G,\varrho}$.

Let $\phi_0=\lambda(\Omega)^{-1}\mathbf1_\Omega$. Then $\phi_0^{(2)}\in
\mathcal S ^{\Omega,K}_{G,\varrho}$ so that $\Phi^{\Omega,K}_{G,\varrho}$
takes finite values. Clearly, $n\mapsto\Phi^{\Omega,K}_{G,\varrho}(n)$
is nonincreasing because $n\mapsto\phi^{(2n)}(e)$
is nonincreasing when $\phi$ is symmetric. By definition,
$\Phi^{\Omega,K}_{G,a\varrho}=\Phi^{\Omega,K}_{G,\varrho}$ for any $a>0$.
A priory, it is possible that $\Phi^{\Omega,K}_{G,\varrho}\equiv0$,
but in
many cases, this possibility can be ruled out so that $\Phi^{\Omega
,K}_{G,\varrho}$
is actually
meaningful and contains information. As indicated below, the choice of
$\Omega$ and $K$
in this definition is mostly irrelevant.

The following proposition contains basic (but not entirely obvious)
properties of $\Phi^{\Omega,K}_{G,\varrho}$
that indicate that Definition~\ref{def1} is quite reasonable. Because
of this
proposition, we will often omit the reference to $\Omega$ and $K$ in
$\Phi^{\Omega,K}_{G,\varrho}$ and write
\[
\Phi^{\Omega,K}_{G,\varrho}=\Phi_{G,\varrho}.
\]

\begin{pro} \label{pro1}
Let $G$ be a locally compact unimodular group. Let $\varrho\dvtx G\ra
[0,+\infty]$ be a measurable function and fix a compact symmetric
neighborhood $\Omega$ of $e$ in $G$
such that $\lambda(\Omega)\ge1$ and $\sup_{\Omega^2}\{\varrho\}>0$.
Fix $K>1$.
\begin{itemize}
\item If there exists a constant $C$ such that, for all $x,y\in G$,
$\varrho(xy)\le C(\varrho(x)+\varrho(y))$ then,
for each integer $n$, $\Phi^{\Omega,K}_{G,\varrho}(n)>0$.
\item For any symmetric continuous
probability density $\phi$ with finite $\varrho$-moment, that is, such that
$\int\varrho\phi \,d\lambda <\infty$, there are a positive constant
$c=c(\phi)$ and an integer $k=k(\phi)$ such that, $\forall n$, $\phi
^{(2n)}(e)\ge c\Phi^{\Omega,K}_{G,\varrho}(kn)$.
\item For $i=1,2$, fix constants $K_i>1$ and compact symmetric
neighborhoods $\Omega_i$ of $e$ in $G$ with $\lambda(\Omega_i)\ge1$.
Let $\varrho_i$, $i=1,2$, be nonnegative measurable functions on $G$
such that $a\varrho_1\le\varrho_2\le A\varrho_1$ for some $a,A\in
(0,\infty)$ and $\sup_{\Omega_i^2}\{\varrho_i\}\in(0,\infty)$. Then,
we have
\[
\Phi^{\Omega_1,K_1}_{G,\varrho_1}\simeq\Phi^{\Omega_2,K_2}_{G,\varrho_2}.
\]
\end{itemize}
\end{pro}

For general $\varrho$, we do not expect to be
able to give a precise bound on $\Phi_{G,\varrho}$, even in the case of
Abelian groups such as $\mathbb Z^d$.

A more reasonable question is to try to understand
$\Phi_{G,\varrho}$ when $\varrho=\rho_G$ and $\rho$ belongs to a
specific family of examples
such as the families
$\rho_\alpha, \rho^{\exp}_{c,\alpha}$, or $\rho^{\log}_\alpha$
mentioned above.
Indeed, in such cases, the function $\Phi_{G,\rho_G}$
(or, perhaps, its equivalence class under the equivalence relation
$\simeq$)
can be thought of as a
group invariant describing the fastest possible decay of
the probability of return of a random walk driven by a
symmetric measure with finite $\rho_G$-moment.
In this restricted context, one may hope to estimate
$\Phi_{G,\rho_G}$ in terms of the function $\Phi_G$ in~(\ref{Phi})
and the function $\rho$. Further, it is an interesting natural question to
ask whether or not all/some of the invariants $\Phi_{G,\rho_G}$ are actually
already determined by $\Phi_G$. This appears to be a rather subtle question.

Another\vspace*{1pt} interesting question raised by Definition~\ref{def1}
is the question of describing classes of measures that are in
$\mathcal S^{\Omega,K}_{G,\varrho}$ and approach\vspace*{1pt}
the extremal behavior described by
$\Phi_{G,\rho_G}$. What is the typical ``shape'' of an almost optimal density?
For instance, should we expect these densities to include
densities that are roughly ``radial'' in terms of the given word-length
\mbox{$|\cdot|$}?
Can we obtain almost extremal densities as convex combinations of
the convolution powers of the uniform probability on a compact symmetric
generating neighborhood of the identity element in~$G$?

Let us observe that determining the exact behavior of $\Phi_{G,\rho_G}$
is a delicate task, even for $G=\mathbb Z$ and $\rho(x)=(1+|x|)^\alpha$,
$\alpha\in(0,2)$. Hence, it is useful and natural to introduce
simplified invariants by comparing $\Phi_{G,\varrho}$
to certain scales of functions. The following definition introduces
a sample of such simplified invariants.
\begin{defin} For $G$ and $\varrho$ as in Definition~\ref{def1}, define:
\begin{longlist}[(3)]
\item[(1)] The power decay invariant,
\[
\mathrm{power}(G,\varrho)=\inf\Bigl\{ \gamma\in(0,\infty)\dvtx \sup_{n}
\{n^\gamma\Phi_{G,\varrho}(n)\}=\infty\Bigr\}.
\]
\item[(2)] The exponential-polylog decay invariant,
\[
\mathrm{exp}\mbox{-}\mathrm{plg}(G,\varrho)= \inf\Bigl\{ \gamma\in(0,\infty)\dvtx\inf_{n}
\bigl\{ \bigl(\log(e+ n)\bigr)^{-\gamma}
\log\bigl(1/\Phi_{G,\varrho}(n)\bigr)\bigr\}=0\Bigr\}.
\]
Computing this quantity is of interest when $\operatorname{power}(G,\varrho
)=\infty$.
\item[(3)] The exponential-power decay invariant,
\[
\mathrm{exp}\mbox{-}\mathrm{pow}(G,\varrho)= \inf\Bigl\{ \gamma\in(0,1]\dvtx\inf_{n}
\bigl\{ n^{-\gamma} \log\bigl(1/\Phi_{G,\varrho}(n)\bigr)\bigr\}=0\Bigr\}.
\]
Again, computing this quantity is of interest when
$\mathrm{exp}\mbox{-}\mathrm{plg}(G,\varrho)=\infty$.
\end{longlist}
\end{defin}
%
\subsection{A sample of illustrative results}
Throughout this subsection we assume that $G$ is compactly generated
and that
$\Phi_G$ is the function given by~(\ref{Phi})
(up to the equivalence relation $\simeq$). With the notation introduced above,
we can state a number of theorems that illustrate the type of results
we obtain in this work. Recall the following notation:
\begin{itemize}
\item$\rho_\alpha(t)= (1+t)^\alpha$, $\alpha\in(0,2)$.
\item$\rho^{\exp}_{c,\alpha}(t)=\exp( c [\log(1+t)]^\alpha)$, $\alpha
\in(0,1)$ and $c>0$;
\item$\rho^{\log}_\alpha(t)= [\log(e+t)]^\alpha$, $\alpha\in(0,\infty)$.
\end{itemize}
\begin{theorem} \label{th-Ipol}
If $G$ has polynomial volume growth of degree $D$, that is, $V(n)\simeq n^D$,
then
\[
\forall\alpha\in(0,2)\qquad \mathrm{power}(G,\rho_\alpha)= D/\alpha
\]
and
\[
\forall\alpha\in(0,1)\qquad \mathrm{exp}\mbox{-}\mathrm{plg}(G,\rho^{\exp}_{c,\alpha})=
1/\alpha.
\]
\end{theorem}

As we shall see, we run into difficulties when estimating
$\mathrm{exp}\mbox{-}\mathrm{pow}(G,\rho^{\log}_\alpha)$.
Assuming $G$ has polynomial volume growth,
we are only able to obtain the estimates
\[
\frac{1}{\alpha+1}\le\mathrm{exp}\mbox{-}\mathrm{pow}(G,\rho^{\log}_\alpha) \le
\frac{1}{\alpha},\qquad \alpha>1.
\]
This indicates that our techniques need to be improved in order to treat
low moment conditions. Indeed, on $\mathbb Z$ (and other Abelian groups),
simple Fourier analysis techniques yield
\[
\mathrm{exp}\mbox{-}\mathrm{pow}(\mathbb Z,\rho^{\log}_\alpha) =
\frac{1}{\alpha+1},\qquad \alpha>0;
\]
see~\cite{BSClmqi}.
\begin{theorem} \label{th-I2}
Assume that the group $G$ has the property that
\[
\forall n\qquad \Phi_G(n)
\ge\exp(-c n^\gamma)
\]
for some $c\in(0,\infty)$ and $\gamma\in(0,1)$.
Then, for any $\alpha\in(0,2)$, there exists $c_1\in(0,\infty)$ such that
\[
\forall n\qquad \Phi_{G,\rho_\alpha}(n)\ge\exp(- c_1 n^{\gamma_\alpha})
\qquad\mbox{where } \gamma_\alpha= \frac{\gamma}{\gamma+(\alpha/2)(1-\gamma)}.
\]
\end{theorem}

So, for instance, for any finitely generated polycyclic group
with exponential volume growth, we have $\gamma=1/3$, and thus
the probability of return of a random walk driven by a symmetric
measure $\mu$ with finite first moment [i.e., $\mu(|\cdot|)<\infty$]
is bounded below by
\[
\mu^{(2n)}(e)\ge\exp(-c_1n^{1/2}).
\]

As indicated by the following results,
the lower bound stated in Theorem~\ref{th-I2} is essentially sharp
in a number of cases.
\begin{theorem}\label{th-I3} Assume that the group $G$ has exponential volume
growth and satisfies
$\forall n, \Phi_G(n)
\ge\exp(-c n^{1/3})$. Then, for each $\alpha\in(0,2)$,
\[
\mathrm{exp}\mbox{-}\mathrm{pow}(G,\rho_\alpha)= \frac{1}{1+\alpha}
\]
and
$\mathrm{exp}\mbox{-}\mathrm{pow}(G,\rho^{\exp}_{c,\beta})=
\mathrm{exp}\mbox{-}\mathrm{pow}(G,\rho^{\log}_\alpha)=1$, $\beta\in(0,1)$, $c>0$,
$\alpha> 2$.
\end{theorem}

Note that the statement that $ \mathrm{exp}\mbox{-}\mathrm{pow}(G,\rho^{\exp}_{c,\beta})=
\mathrm{exp}\mbox{-}\mathrm{pow}(G,\rho^{\log}_\alpha)=1$ for the groups considered
in Theorem~\ref{th-I3} is crude. More detailed results are described in the
core of the paper. For instance, $\mathrm{exp}\mbox{-}\mathrm{pow}(G,\rho^{\exp}_{c,\beta})=1$
can be refined to the much more informative statement that, for any
fixed $c>0$
and $\beta\in(0,1)$, there exist $c_1,c_2\in(0,\infty)$ such that
\[
-n \exp(- c_1(\log n)^\beta) \le
\log\Phi_{G,\rho^{\exp}_{c,\beta}}(n)\le- n \exp(- c_2(\log n)^\beta)
\]
for all $n$ large enough.

The case of the lamplighter groups
$(\mathbb Z/2\mathbb Z)\wr\mathbb Z^d$, the simplest wreath products,
is particularly interesting.
\begin{theorem}\label{th-I4} For $G_d=(\mathbb Z/2\mathbb Z)\wr\mathbb Z^d$,
$d=1,2,\ldots,$ and for $\alpha\in(0,2)$,
\[
\mathrm{exp}\mbox{-}\mathrm{pow}(G_d,\rho_\alpha)= \frac{d}{d+\alpha}.
\]
\end{theorem}
\begin{pf}
The upper bound follows from Theorem
\ref{th-I2}. The lower bound requires an ad hoc argument explained in Section
\ref{sec-Wreath}.
\end{pf}

For the next result, recall that a group $G$ is meta-Abelian if it
contains a
normal Abelian subgroup $A$\vadjust{\goodbreak} such that $G/A$ is Abelian. From the view point
of group theory,
meta-Abelian groups are only ``one step'' removed from being Abelian.
\begin{theorem} Let $G$ be a finitely generated meta-Abelian group. Then
either $G$ has polynomial volume growth and there is an integer $D$
such that
\[
\forall\alpha\in(0,2)\qquad \mathrm{power}(G,\rho_\alpha)= D/\alpha
\]
or there exists an integer $d$ such that
\[
\forall\alpha\in(0,2)\qquad
\frac{1}{1+\alpha}\le\mathrm{exp}\mbox{-}\mathrm{pow}(G,\rho_\alpha)
\le\frac{d}{d+\alpha}.
\]
\end{theorem}
\begin{pf}
Being solvable, finitely generated
meta-Abelian groups either
have polynomial volume growth or exponential volume growth; see \cite
{Milnor,Wolf}.
In the polynomial volume growth case, apply Theorem~\ref{th-Ipol}. For
any group with exponential volume growth,
Theorem~\ref{th-Vol2} gives the lower bound $\mathrm{exp}\mbox{-}\mathrm{pow}(G,\rho
_\alpha)
\ge\frac{1}{1+\alpha}$. By~\cite{PSCwp}, any\vspace*{1pt} meta-Abelian group has
$\Phi_G(n)\ge\break\exp(-Cn^{d/d+2})$ for some integer $d\ge1$.
Thus the upper bound $\mathrm{exp}\mbox{-}\mathrm{pow}(G,\break\rho_\alpha)
\le\frac{d}{d+\alpha}$ follows from Theorem~\ref{th-I2}.
\end{pf}

\subsection{Methodology} We close this introduction by describing in
broad terms
the techniques we will use to prove the results described above. For
the purpose
of this discussion, we focus on the problem of estimating the rate of decay
of convolution powers of symmetric measures having a continuous density
and a finite $\rho_\alpha$-moment,
$\alpha\in(0,2)$ [$\rho_\alpha(x)=(1+|x|)^\alpha$]. We start with a
quick review of classical results
in the context of the lattice $\mathbb Z^d$. In this context,
the literature focuses on local limit theorems, that is,
results that describe the precise asymptotic behavior of $\phi^{(2n)}(x)$.
For instance, if $\phi$ is a symmetric probability density
which has generating support and finite
second moment, $\phi^{(2n)}(0)\sim c(d,\mu) n^{-d/2}$ (e.g.,
\cite{Spitzer}, P9, Section 7). For $\alpha\in(0,2)$,
the simple condition of having a finite $\rho_\alpha$-moment is not sufficient
for the validity of a local limit theorem, even on $\mathbb Z$.

For a symmetric probability density $\phi$ on $\mathbb Z$,
set $G(k)= \sum_{|i|\ge k}\phi(k)$,
$H(k)= k^{-2}\sum_{|i|\le k}i^2\phi(i)$. Then $\phi$ is
in the domain of attraction of a symmetric stable law of index $\alpha$ if
$\lim_{\infty}H/G= \alpha/(2-\alpha)$.
In such a case, a local limit theorem holds
stating that $\phi^{(2n)}(e)\sim c(\alpha,\mu) a_n$
with $a_n$ defined by $Q(a_n)=1/n$, $Q=G+H$; see, for example,
\cite{Fel2,FelB5,Grif,GJP}.
All classical discussions of such results make heavy use of
Fourier transform techniques. It\vspace*{1pt} is easy to use these techniques to see
that if
a symmetric probability density $\phi$ with generating support on
$\mathbb Z^d$
has finite $\rho_\alpha$-moment for some $\alpha\in(0,2)$, then we
must have
$\phi^{(2n)}(0)\ge c(d,\mu) n^{-d/\alpha}$ and
\[
\Phi_{\mathbb Z^d,\rho_\alpha}(n)\ge c(d,\alpha) n^{-d/\alpha}.
\]
As laws that are in the domain of attraction of a symmetric stable law of
index $\beta>\alpha$ have finite $\rho_\alpha$-moment, we also get that
\[
\forall\beta>\alpha\qquad
\Phi_{\mathbb Z^d,\rho_\alpha}(n)\le c_\beta n^{-d/\beta}.
\]
Hence $\operatorname{power}(\mathbb Z^d,\rho_\alpha)=d/\alpha$. Note that determining
the exact behavior of $n\mapsto\Phi_{\mathbb Z,\rho_\alpha}(n)$ appears
to be a somewhat subtle problem and will not be discussed here.

Both the Fourier transform and explicit examples such as symmetric
stable laws are not available on most noncommutative groups so that the
arguments outlined above must be replaced by different ideas. Our
approach is as follows:

\begin{longlist}[(2)]
\item[(1)]
Our lower bounds on $\Phi_{G,\rho_\alpha}$ are obtained and expressed
via the function $\Phi_G$ given by~(\ref{Phi}).
This function $\Phi_G$ describes the decay of convolution powers
of symmetric, nondegenerate
densities with finite second moment on the group $G$.
To transfer the information contained in this
function $\Phi_G$ and make it relevant to the study of the convolution powers
of measures with finite $\rho_\alpha$-moment, we will use a sort of
interpolation argument, the comparison of Dirichlet forms
and the notion of von Neumann trace.
Each one of these ingredients plays a crucial role in obtaining
our lower bounds.

Section~\ref{sec-CompD} contains the
proof of Proposition~\ref{pro1} as well as an interesting
and important variation on Definition~\ref{def1}. It also
develops the key interpolation argument which leads to the
comparison of important quadratic forms including the Dirichlet forms
of the probability measures we want to study.

The role of the notion of von Neumann trace is explained in the
\hyperref[app]{Appendix} where related needed material is described. The results
developed in the \hyperref[app]{Appendix} are the tools that allow us to turn
the comparison of quadratic forms obtained in Section~\ref{sec-CompD}
into lower bound for $\Phi_{G,\rho_\alpha}$.

\item[(2)] To obtain upper bounds on $\Phi_{G,\rho_\alpha}$, it suffices
to exhibit some probability densities satisfying the desired moment condition
and whose convolution powers can be estimated. On a general noncommutative
group, this is not necessarily an easy task.
One possible technique---discrete subordination---uses
Bernstein functions
to produce probability densities on $G$ that include laws
that can be thought of as analogs of symmetric stable laws.
The decay of the convolution powers
of these laws can be precisely expressed and controlled in terms of the group
invariant $\Phi_G$ at~(\ref{Phi}), and this technique is quite
interesting in
its own right. This idea, which the authors developed for the purpose
of the
present paper, is presented in detail in~\cite{BSCsubord}.
We will use some of the results of~\cite{BSCsubord}.
However, the moment properties of these subordinated laws are directly
related to the rate of escape of the basic simple random walk on the
underlying group. In particular, for groups with a rate of escape
that is faster than the classical $\sqrt{n}$, the moment conditions
satisfied by these subordinated laws are not what one would expect
from a simplistic analogy with the classical case of $\mathbb Z$.
For instance, on a finitely generated group with linear rate of escape,
the ``symmetric stable law''
of exponent $\beta\in(0,2)$ [by definition, the law obtained via discrete
$(\beta/2)$-subordination from the
law of simple random walk] will only have a finite $\rho_{\alpha}$-moment
for $\beta>2\alpha$ (instead of $\beta>\alpha$ in the classical case).
What this means is that, in general,
upper bounds obtained by using~\cite{BSCsubord}
will not match closely the lower bounds
discussed in (1) above. They will only do so if there exists
a simple random walk on $G$ that
has a rate of escape of type $\sqrt{n}$ as in the classical case of
$\mathbb Z$.
This is a subtle requirement since it is not known whether or not
all random walks associated with finite symmetric generating sets on a given
finitely generated group have the same rate of escape.

\item[(3)] There is a more elementary way to produce probability
distributions with
finite $\rho_\alpha$-moment and whose convolution powers can be estimated.
This technique goes back to~\cite{SCconv,Varconv}.
It is revisited in Section~\ref{sec-UpVol}. It works well for groups where
the invariant $\Phi_G$ behaves precisely as predicted by the available
upper bounds based on volume growth (e.g., polycyclic groups).
It does not work well for
wreath products such as $(\mathbb Z/2\mathbb Z)\wr\mathbb Z^d$, $d\ge2$.
For such groups neither of the techniques in (2) or (3) produce upper bounds
on $\Phi_{G,\rho_\alpha}$ matching the lower bounds obtained via (1).
Nevertheless, Section~\ref{sec-Wreath} shows that the lower bounds obtained
via (1) are essentially tight even in the case of these wreath products.
This requires an ad hoc argument that takes advantage of the precise
structure of these groups.
\end{longlist}

\section{Comparisons of Dirichlet forms}
\label{sec-CompD}

This section develops the key technique that we use to obtain
lower bounds on the functions $\Phi_{G,\varrho}$ introduced in
Definition~\ref{def1},
namely, comparison of Dirichlet forms. The first
subsection contains
simple results that show that the object introduced in Definition~\ref{def1},
$\Phi_{G,\varrho}$, has some nice stability properties. The second subsection
develops a somewhat sophisticated comparison between certain quadratic forms.
It plays a central role in our results.

It is useful to introduce the following somewhat subtle modification
of Definition~\ref{def1} in which a ``weak moment condition,''
$W(\varrho,\mu)<\infty$, replaces the
``strong moment condition'' $\mu(\varrho)<\infty$. For any probability measure
$\mu$ and $\varrho\dvtx G\ra[0,\infty)$, $W(\varrho, \mu)$ is defined by
\[
W(\varrho,\mu)= \sup_{s>0} \{s \mu( \varrho>s)\}.
\]

\begin{defin} \label{def2}
Let $G$ be a locally compact unimodular group.
Fix a measurable function $\varrho\dvtx G\ra[0,+\infty]$.
Fix a compact symmetric neighborhood $\Omega$ of $e$ in $G$
such that $\lambda(\Omega)\ge1$ and $\sup_{\Omega^2}\{\varrho\}>0$.
For $K> 1$, let $\widetilde{\mathcal S}^{\Omega,K}_{G,\varrho}$
be the set of all symmetric continuous
probability densities $\phi$ on $G$ with the properties that
$\|\phi\|_\infty\le K$ and
$W(\varrho,\phi \,d\lambda) \le K \sup_{\Omega^2}\{\varrho\}$. Set
\[
\widetilde{\Phi}^{\Omega,K}_{G,\varrho}\dvtx
n\mapsto\widetilde{\Phi}^{\Omega,K}_{G,\varrho}(n):=
\inf\bigl\{
\phi^{(2n)}(e) \dvtx
\phi\in\widetilde{\mathcal S}^{\Omega,K}_{G,\varrho}\bigr\}.
\]
\end{defin}

Obviously,
%
\begin{equation}
\Phi^{\Omega,K}_{G,\varrho}\ge \widetilde{\Phi}^{\Omega,K}_{G,\varrho}.
\end{equation}
In the classical case of $\mathbb R^d$ or $\mathbb Z^d$ with
$\varrho=\rho_\alpha(|\cdot|)=(1+|\cdot|)^\alpha$,
as long as $\alpha\in(0,2)$, we have
\[
\widetilde{\Phi}^{\Omega,K}_{\mathbb Z^d,\rho_\alpha}(n)
\simeq n^{-d/\alpha},
\]
whereas it is not easy to estimate
$\Phi^{\Omega,K}_{\mathbb Z^d,\rho_\alpha}(n)$ precisely
(see the comments made in the \hyperref[sec-Intr]{Introduction}).
Interestingly enough, for $\alpha=2$, we have (see~\cite{Grif})
\[
\Phi^{\Omega,K}_{\mathbb Z^d,\rho_2}(n)\simeq n^{-d/2},\qquad
\widetilde{\Phi}^{\Omega,K}_{\mathbb Z^d,\rho_2}(n)\simeq(n\log n)^{-d/2}
.
\]

\subsection{\texorpdfstring{Some basic stability results for $\Phi_{G,\varrho}$}
{Some basic stability results for Phi G, rho}}

By definition, a continuous symmetric probability density $\phi$ such that
\[
\|\phi\|_\infty\le K \quad\mbox{and}\quad
\int\varrho\phi \,d\lambda\leq K\sup_{\Omega^2}\{\varrho\}
\]
must satisfy
\[
\phi^{(2n)}(e)\ge \Phi^{\Omega,K}_{G,\varrho}(n).
\]
It is natural to ask what can be
said of a symmetric probability density $\phi\in L^2(G)$ such that
$\int\varrho\phi \,d\lambda<\infty$. This section gives a reassuring
answer to this question and proves the results stated in Proposition
\ref{pro1}.
We need the following elementary fact.
\begin{pro}
Let $G$ be a locally compact unimodular group. Assume that
$\varrho\dvtx G\ra[0,\infty]$ is a measurable function with the property that
there exists $C\in[1,\infty)$ such that
\[
\forall x,y\in G\qquad \varrho(xy)\le C\bigl(\varrho(x)+\varrho(y)\bigr).
\]
If $\mu$ is a probability measure satisfying $\mu(\varrho)<\infty$, then
\[
\mu^{(n)}(\varrho) \le nC^{n-1}\mu(\varrho),\qquad n=1,2,\ldots.
\]
Further, we have
\[
W\bigl(\varrho,\mu^{(n)}\bigr) \le
n (2C)^{n-1} W(\varrho,\mu).
\]
\end{pro}
\begin{pf} By definition of the convolution product,
for any two measures $\mu,\nu$,
\[
\mu*\nu(\varrho)= \int_{G\times G}\varrho(xy)\, d\mu(x)\,d\nu(y).
\]
If $\mu,\nu$ are probability measures, since $\varrho(xy)\le C(\varrho
(x)+\varrho(y))$, we obtain
\[
\mu*\nu(\varrho) \le C \bigl(\mu(\varrho)+\nu(\varrho)\bigr).\vadjust{\goodbreak}
\]
The\vspace*{1pt} inequality $\mu^{(n)}(\varrho) \le nC^{n-1}\mu(\varrho)$
follows by induction. To obtain the inequality regarding
$W(\varrho,\mu^{(n)})$ observe that
\[
\{(x,y)\dvtx \varrho(x,y)>s\}\subset\{(x,y)\dvtx \varrho(x)>s/(2C)\}
\cup\{(x,y)\dvtx\rho(y)>s/(2C)\}.
\]
Hence, for any two probability measures $\mu,\nu$, we have
\begin{eqnarray*}
\mu*\nu(\{\varrho>s\}) &=&
\int_{\{(x,y)\dvtx \varrho(xy)>s\}} d\mu(x)\,d\nu(y)\\
&\le& \mu(\{\varrho>s/2C\})+\nu\bigl(\{\varrho>s/(2C)\}\bigr).
\end{eqnarray*}
This yields $W(\varrho,\mu*\nu)\le2C(W(\varrho,\mu)+W(\varrho,\nu))$
and the desired result follows by induction.
\end{pf}
\begin{cor} Let $\varrho,\Omega,K$ be as in Definition~\ref{def1}.
Assume that $\varrho$
tends to infinity at infinity and satisfies
\[
\forall x,y\in G\qquad
\varrho(xy)\le C\bigl(\varrho(x)+\varrho(y)\bigr).
\]
Then $\Phi^{\Omega,K}_{G,\varrho}(n) >0$ and $\widetilde{\Phi}^{\Omega
,K}_{G,\varrho}(n) >0$,
for all $\Omega,K,n$.
\end{cor}
\begin{pf} We prove the result for $\Phi^{\Omega,K}_{G,\varrho}$ (the
case of $\widetilde{\Phi}^{\Omega,K}_{G,\varrho}$ is similar).
Let $\varrho_0=\sup_{\Omega^2}\{\varrho\}$.
Let $\phi$ be a symmetric continuous probability density
in $\mathcal S^{\Omega,K}_{G,\varrho}$, that is,
such that $\|\phi\|_\infty\le K$ and
$\int\phi\varrho \,d\lambda\le K\varrho_0$.
Then $\int\phi^{(2n)}\varrho \,d\lambda\le2KnC^{2n-1}\varrho_0$.
Further, for any $N$,
\[
\int_{\varrho>N} \phi^{(2n)}\,d\lambda\le
N^{-1}\int\phi^{(2n)}\varrho \,d\lambda\le \frac{2KnC^{2n-1}\varrho_0}{N}.
\]
Since $\phi^{(2n)}$ attains its maximum at $e$, we obtain that
\[
\phi^{(2n)}(e)\ge\frac{1}{\lambda(\varrho\le N)}
\int_{\varrho\le N} \phi^{(2n)}\,d\lambda\ge
\frac{1- {2KnC^{2n-1}\varrho_0}/{N}}{\lambda(\varrho\le N)}.
\]
As $\varrho$ tends to infinity at infinity, $\lambda(\varrho\le N)$ is finite
for any finite $N$. Hence, for $N= 4KnC^{2n-1}\varrho_0$, we obtain
a uniform positive lower bound on
$\phi^{(2n)}(e)$ for all $\phi\in\mathcal S^{\Omega,K}_{G,\varrho}$.
\end{pf}
\begin{pro}
Let $G$ be a locally compact unimodular group.
Let $\varrho,\Omega,K$ be as in Definition~\ref{def1}.
Let $\phi$ be a symmetric continuous probability density.
\begin{itemize}
\item Assume that
$\int\varrho\phi \,d\lambda<\infty$. Then there exist
$c_1=c_1(\phi)>0,c_2=c_2(\phi)\in\mathbb N$ such that
\[
\forall n=1,2,\ldots\qquad
\phi^{(2n)}(e)\ge c_1\Phi^{\Omega,K}_{G,\varrho}(c_2n).
\]
\item Assume that
$W(\varrho,\phi \,d\lambda)<\infty$.
Then there exist
$c_1=c_1(\phi)>0,c_2=c_2(\phi)\in\mathbb N$ such that
\[
\forall n=1,2,\ldots\qquad
\phi^{(2n)}(e)\ge c_1\widetilde{\Phi}^{\Omega,K}_{G,\varrho}(c_2n).\vadjust{\goodbreak}
\]
\end{itemize}
\end{pro}
\begin{pf} The proofs of the two statements are similar and we only
give the proof under the condition
$\int\varrho\phi \,d\lambda<\infty$. Let $\phi_0=\lambda(\Omega
)^{-1}\mathbf1_\Omega$.
Obviously, since $\lambda(\Omega)\ge1$, we have $\|\phi_0^{(2)}\|
_\infty\le1$ and
$\int\varrho\phi^{(2)}_0\,d\lambda\le\sup_{\Omega^2}\{\varrho\}$. By hypothesis,
\[
M=\max\biggl\{\|\phi\|_\infty,\Bigl(\sup_{\Omega^2}\{\varrho\}\Bigr)^{-1}\int\varrho
\phi
\,d\lambda\biggr\}<+\infty.
\]
If $M\le K$, the result is clear. If not then $M>K>1$. In this case,
set $\alpha=(K-1)/(M-1)\in(0,1)$, and observe that
the symmetric continuous
probability density
$\phi_1=\alpha\phi+(1-\alpha)\phi^{(2)}_0$ satisfies
$\|\phi_1\|_\infty\le K$ and $\int\varrho\phi_1\,d\lambda\le K\sup
_{\Omega^2}\{\varrho\}$. Thus,
\[
\forall n\qquad \phi_1^{(2n)}(e)\ge\Phi^{\Omega,K}_{G,\varrho}(n).
\]
Further, by construction, the Dirichlet forms
$\mathcal E=\mathcal E_{\phi \,d\lambda} $ and
$\mathcal E_1=\mathcal E_{\phi_1 \,d\lambda}$ [see~(\ref{Emu})] satisfy
$\mathcal E\le(1/\alpha) \mathcal E_1$.
In terms of the convolution operator $R_\phi$ (convolution on the right
by $\phi$) acting on $L^2(G)$, this is equivalent
to say that
\[
I-R_\phi\le(1/\alpha)( I-R_{\phi_1}).
\]
By Corollary~\ref{cor-comp1}, this implies that
$\phi^{(2n)}(e)\ge c_1 \phi_1^{(2 c_2 n )}(e)$,
for some $c_1>0$ and $c_2\in\mathbb N$.
\end{pf}
\begin{rem} If, for all $x,y\in G$, $\varrho(xy)\le C(\varrho(x)+\varrho(y))$,
then any symmetric probability density $\phi\in L^2(G)$
(not necessarily continuous) with finite $\varrho$-moment satisfy
\[
\phi^{(2n)}(e)\ge c_1\Phi^{\Omega,K}_{G,\varrho}(c_2 n)
\]
for some $c_1=c_1(\phi)>0$ and $c_2=c_2(\phi)\in\mathbb N$.
Indeed, it suffices to apply the previous result to $\phi*\phi$ which
is continuous and
also has finite $\varrho$-moment.
\end{rem}
\begin{pro}
Let $G$ be a locally compact unimodular group.
For $i=1$, $2$, fix constants $K_i>1$ and compact neighborhoods $\Omega_i$
of $e$ in $G$ with $\lambda(\Omega_i)\ge1$. Let
$\varrho_i\dvtx G\ra[0,\infty)$, $i=1,2$, be measurable functions
such that $a\varrho_1\le\varrho_2\le A\varrho_1$ for some $a,A\in
(0,\infty)$ and $\sup_{\Omega_i^2}\{\varrho_i\}\in(0,\infty)$. Then,
we have
\[
\Phi^{\Omega_1,K_1}_{G,\varrho_1}\simeq\Phi^{\Omega_2,K_2}_{G,\varrho_2},\qquad
\widetilde{\Phi}^{\Omega_1,K_1}_{G,\varrho_1}\simeq\widetilde{\Phi
}^{\Omega_2,K_2}_{G,\varrho_2}.
\]
\end{pro}
\begin{pf}We treat the case of the function $\Phi$. The case of
$\widetilde{\Phi}$ is similar.
Set $M_i=\sup_{\Omega_i^2}\{\varrho_i\}\in(0,\infty)$.
Let\vspace*{-1pt} $\phi_1\in\mathcal S^{\Omega_1,K_1}_{G,\varrho_1}$.
Let $\phi_0=\lambda(\Omega_2)^{-1}\mathbf1_{\Omega_2}$ and set
$\phi_2= \alpha\phi_1+ (1-\alpha)\phi_0^{(2)}$,
for some $\alpha\in(0,1] $ to be chosen later.
This continuous symmetric probability density satisfies
\[
\|\phi_2\|_\infty\le\alpha K_1+ (1-\alpha) \quad\mbox{and}\quad
\int\phi_2 \varrho_2 \,d\lambda\le\alpha A K_1M_1 + (1-\alpha)M_2.
\]
It follows that, for $\alpha$ close enough to $1$, we have $\phi_2\in
\mathcal S^{\Omega_2,K_2}_{G,\varrho_2}$.
Indeed, picking $\alpha=\min\{1,(K_2-1)/(K_1-1),
(K_2-1)M_2/AK_1|M_2-M_1| \}$ will work.

As in the previous proof, setting $\mathcal E_i=
\mathcal E_{\phi_i \,d\lambda}$ [see~(\ref{Emu})], we find that
$\mathcal E_1\le(1/\alpha)\mathcal E_2$.
By Corollary~\ref{cor-comp1}, this implies that
there exists $c>0$ and an integer $k$ such that
\[
\phi_1^{(2n)}(e)\ge c \phi_2^{(2kn)}(e).
\]
Further $c$ and $k$ depend only $K_1,K_2, M_1,M_2$ and $A$. Hence
\[
\forall n\qquad \Phi^{\Omega_1,K_1}_{G,\varrho_1}(n)\ge c \Phi^{\Omega
_2,K_2}_{G,\varrho_2}(kn)
\]
as desired. Using the symmetry of the hypotheses, the reverse
inequality holds as well.
\end{pf}

\subsection{An abstract interpolation/comparison result}
The results developed in this key section make use of a given nonnegative
self-adjoint operator $(A,D_A)$ on $L^2(G)$ with associated semigroup
$H_t=e^{-tA}$, $t\ge0$, which is assumed to be, in some sense, well understood.
In applications, $H_t$ will actually be a symmetric Markov semigroup, and
$\|A^{1/2} f\|_2^2$ will be a Dirichlet form. We assume that $A$
(and thus also $H_t$) commutes with left translations in $G$. Namely,
for $f\in L^2(G)$ and $h\in G$, set $\tau_h f=f(h \cdot)\in L^2(G)$.
We assume that
$A$ has the property that $ f\in D_A$ implies
$\tau_h f\in D_A$ and $A (\tau_h f)= \tau_h(Af)$ for any $h\in G$.

As mentioned above, we think of the semigroup $H_t$ as a basic object
which is well understood. The key idea is that we then also understand quite
well the semigroups generated by certain functions $\psi(A)$ of $A$.
The class of functions $\psi$ of interest to us here
is the class of those functions that admit the Laplace-type representation
%
\begin{equation}
\label{omegapsi}
\psi(\lambda)=\lambda^2\int_0^\infty e^{-\lambda s} \omega(s)\,ds
\qquad\mbox{with } \omega\ge0.
\end{equation}
The simplest example of such function is $\psi\dvtx\lambda\mapsto\lambda
^\alpha$,
$\alpha\in(0,1)$, which is obtained by picking
$\omega(s)= c_\alpha s^{1-\alpha}$, $c_\alpha=1/\Gamma(2-\alpha)$.
By spectral theory the $L^2(G)$-domain of $\psi(A)^{1/2}$
is the set of functions
$f\in L^2(G)$ such that
%
\begin{equation}\label{psiomega}
\|\psi(A)^{1/2}f\|_2^2= \int_0^\infty\|AH_{s/2}f\|_2^2 \omega
(s)\,ds<\infty.
\end{equation}
It is easy to see that $\|\psi(A)^{1/2}f\|^2_2<\infty$ whenever
$f\in D_A$ and $\omega(s)\le C(1+s)$ (in fact, $f\in D_{A^{1/2}}$ suffices).
It follows that $\psi(A)^{1/2}$ is densely defined and self-adjoint whenever
$\omega(s)\le C(1+s)$.

Next, we introduce a key assumption about $(A,D_A)$. This assumption is
expressed in term of a given positive (measurable) function
$\delta\dvtx G\mapsto[0,\infty)$.
It captures a fundamental relation between the $L^2$-variation of $f$ and
$\|A^{1/2}f\|_2$.
Namely, setting
\[
f_h(x)=f(xh), \qquad f\in L^2(G), x,h\in G,\vadjust{\goodbreak}
\]
we assume that there exists a constant $C_0\in[1,\infty)$ such that
%
\begin{equation}\label{MVrho}
\forall f\in D_A, \forall h\in G\qquad
\biggl(\int_G |f_h-f|^2\,d\lambda\biggr)^{1/2}\le C_0\delta(h) \|A^{1/2}f\|_2.
\end{equation}

Finally, for any probability measure $\mu$ on $G$, we set
%
\begin{equation}\label{Emu}\quad
\forall f\in L^2(G)\qquad \mathcal E_\mu(f,f)=\frac{1}{2}\int_G\int_G
|f(xy)-f(x)|^2 \,d\lambda(x)\,d\mu(y).
\end{equation}
When $\mu$ is symmetric, $\mathcal E_\mu$ is the
Dirichlet form associated with $\mu$ and
$\mathcal E_\mu(f,\break f)=\langle f-f*\mu,f\rangle$.
\begin{theorem}\label{th-interpol}
Referring to the setting and notation introduced above,
consider a pair of nonnegative increasing functions $\omega$, $\psi$
related by
(\ref{omegapsi}).
Assume that
$s\mapsto\omega(s)/s$ is decreasing, and set
%
\begin{equation}
\xi(t)= \int_0^t
\biggl(\frac{s}{\omega(s)}\biggr)^{1/2}\,\frac{ds}{s},\qquad
\zeta(t)=t^{1/2}\int_t^\infty\frac{ds}{s\omega(s)^{1/2}}.
\end{equation}

Let $\rho\dvtx[0,\infty)\ra[1,\infty)$ be an increasing function such that,
for all $t\ge0$,
%
\begin{equation} \label{rho1}
\frac{t\max\{\xi(t^2),\zeta(t^2)\}}{\omega(t^2)^{1/2}}\le C^2_1\rho(t).
\end{equation}
%

Then, if $A$ satisfies~(\ref{MVrho}) and $\mu$ is such that
$\mu(\rho\circ\delta)<\infty$, we have
\[
\mathcal E_\mu(f,f) \le8C_0^2C_1
\mu(\rho\circ\delta) \|\psi(A)^{1/2}f\|_2^2,\qquad f\in D_A.
\]
\end{theorem}
\begin{rem} When $\mu$ is symmetric, $\mathcal E_\mu$ is a Dirichlet form.
In general, $f\mapsto\|\psi(A)^{1/2}f\|_2^2$ is not a Dirichlet form.
If we assume that $-A$ is the infinitesimal generator of
a symmetric Markov semigroup, then $f\mapsto\|\psi(A)^{1/2}f\|_2^2$ is
a Dirichlet form if we assume that $\psi$ is a Bernstein function;
see~\cite{BSCsubord,NJ}. This will not play an important role in this paper
but~\cite{BSCsubord}, Theorem 2.5, shows that it is often possible to
choose $\psi$ to be a Bernstein function.
\end{rem}
\begin{rem}\label{rem-interpol} The functions $\xi, \zeta$ are always
greater or equal to $(t/\omega(t))^{1/2}$.
The typical functions $\omega$ of interest to us are such that
$\omega(s)\ge\eta s^{1-\varepsilon}$ in $(0,1)$ with $\varepsilon\in(0,1)$
and $\omega(s)\simeq s/\beta(s)$
at infinity with $\beta$ an increasing regularly varying function of index
in $[0,1)$. If $\beta$ has index in $(0,1)$, then
\[
\xi(t)\simeq\zeta(t) \simeq
\biggl(\frac{t}{\omega(t)}\biggr)^{1/2} \qquad\mbox{at infinity}
\]
and~(\ref{rho1}) can be replaced by
\[
\frac{t^2}{\omega(t^2)} \le C_1^2\rho(t).
\]
If, instead, $\beta$ is slowly varying then it is still the case that
$\zeta(t)\simeq(t/\omega(t))^{1/2}$ at infinity but, for $t$ large enough,
\[
\bigl(t/\omega(t)\bigr)^{1/2}\ll
\xi(t)\le C \log(e+t) \bigl(t/\omega(t)\bigr)^{1/2}.
\]
In this case, $\max\{\zeta(t^2),\xi(t^2)\}=\xi(t^2)$ for $t$ large enough.

If $\beta$ has index $1$ and is of the form $\beta(t)= t/\ell(t)$
with $\ell(t)$ slowly varying at infinity then $\omega(t)\simeq\ell(t)$,
$\xi(t)\simeq(t/\omega(t))^{1/2}$ but $ (t/\omega(t))^{1/2}\ll \zeta
(t) $
and, in fact, $\zeta(t)$ might be infinite unless further assumptions
are made on $\ell$.
\end{rem}
\begin{pf*}{Proof of Theorem~\ref{th-interpol}}
Let $f\in D_A$ and write
\begin{eqnarray*}
g(h)&=&\|f_h-f\|_2,\\
f_h-f&=& ([H_tf]_h-H_tf) + ([f-H_tf]_h) -(f-H_tf).
\end{eqnarray*}
Since $f-H_tf=\int_0^t AH_s f \,ds$, we have
%
\begin{equation}\label{p1-1}
\|([f-H_tf]_h) -(f-H_tf)\|_2\le2 \int_0^t \|AH_sf\|_2 \,ds.
\end{equation}
Using~(\ref{MVrho}), we also have
%
\begin{equation}\label{p1-2}
\|[H_tf]_h-H_tf\|_2\le C_0 \delta(h)\|A^{1/2}H_tf\|_2.
\end{equation}
Further,
%
\begin{equation}\label{p1-3}\qquad
\|A^{1/2}H_tf\|_2= \biggl\| \int_t^\infty A^{1/2} AH_sf\,ds\biggr\|_2
\le\int_{t}^\infty(es)^{-1/2}\|AH_{s/2} f\|_2\,ds.
\end{equation}
Here we have used the inequalities
\[
\| A^{1/2}H_s f\|_2\le\|A^{1/2}H_{s/2}\|_{2\ra2}
\|H_{s/2}f\|_2
\]
and (by spectral theory)
\[
\|A^{1/2}H_{s}\|_{2\ra2}\le\max_{a>0}\{a^{1/2}e^{-sa}\}= (2es)^{-1/2}.
\]
Putting together inequalities~(\ref{p1-1}),~(\ref{p1-2}) and~(\ref{p1-3})
yields
\[
g(h)\le2\int_0^t\|AH_sf\|_2 \,ds+ C_0\delta(h)\int_t^\infty(es)^{-1/2}\|
AH_{s/2}f\|_2\,ds.
\]
Pick $t= \tau(h)=\max\{1,\delta(h)^2\}$, set $\theta=\max\{\xi,\zeta\}$
and write
\[
\frac{g(h)}{\theta\circ\tau(h)}\le2C_0
\int_0^\infty K(h,s)([s\omega(s)]^{1/2}
\|AH_{s/2} f\|_2) \,\frac{ds}{s},
\]
where $K$ is the kernel on $G\times(0,\infty)$ given by
\[
K(h,s)= \frac{s^{1/2}}{\theta\circ\tau(h)\omega(s)^{1/2}}
\bigl(\mathbf1_{(0,\tau(h))}(s)+ \delta(h)s^{-1/2}\mathbf1_{[\tau(h),
\infty)}(s)\bigr).
\]
Consider this kernel as defining an integral operator
\begin{eqnarray*}
K\dvtx L^2\biggl((0,\infty),\frac{ds}{s}\biggr)
&\ra& L^2(G,[\theta\circ\tau]^2 \,d\mu),
\qquad u\mapsto Ku,
\\
Ku(h)&=&\int_0^\infty K(h,s)u(s)\,\frac{ds}{s}.
\end{eqnarray*}
Assuming that this operator is bounded with norm $N_*$, we obtain
%
\begin{eqnarray} \label{Epsi}
\mathcal E_\mu(f,f)&=&\int_G|g|^2 \,d\mu\nonumber\\
&\le& 4C_0^2N_* ^2\int_0^\infty\|AH_{s/2}f\|_2^2\omega(s) \,ds\\
&=&
4C_0^2N^2_*\|\psi(A)^{1/2}f\|_2^2.\nonumber
\end{eqnarray}
A standard interpolation argument gives
\[
N_*^2\le\biggl( \sup_{h\in G}\int_0^\infty K(h,s)\,\frac{ds}{s} \biggr)
\biggl( \sup_{s>0} \int_G K(\cdot,s) [\theta\circ\delta]^2 \,d\mu \biggr)
\]
and we have
\begin{eqnarray*}
\int_0^\infty K(h,s)\,\frac{ds}{s}&=&
\frac{1}{\theta(\tau(h))}\int_0^{\tau(h)}\,\frac{ds}{[s\omega(s)]^{1/2}}
+ \frac{\delta(h)}{\theta(\tau(h))}\int_{\tau(h)}^\infty\frac
{ds}{s\omega(s)^{1/2}},
\\
\int_GK(\cdot,s)[\theta\circ\tau]^2 \,d\mu&=&
\frac{s^{1/2}}{\omega(s)^{1/2}}\int_{\{\tau> s\}}
[\theta\circ\tau] \,d\mu+
\frac{1}{\omega(s)^{1/2}}\int_{\{\tau\le s\}}
\delta[\theta\circ\tau] \,d\mu.
\end{eqnarray*}
By the definitions of $ \xi,\zeta$ and $\theta$,
\[
\sup_{h\in G}\biggl\{\int_0^\infty K(h,s)\,\frac{ds}{s}\biggr\}\le2.
\]
Further, since we assume that
$s\mapsto\omega(s)$  is increasing and
$s\mapsto\omega(s)/s$  decreasing,~(\ref{rho1}) yields
\[
\sup_{s>0}\biggl\{\int_GK(\cdot,s) [\theta\circ\tau]^2 \,d\mu\biggr\}
\le C_1
\int\rho\circ\delta \,d\mu.
\]
This gives the desired result.
\end{pf*}

This proof admits the following result as a corollary.
\begin{theorem}\label{th-interpol2}
Referring to the setting and notation introduced above,
consider a pair of smooth nonnegative increasing functions
$\omega$, $\psi$ related by~(\ref{omegapsi}).
Fix $\alpha\in(0,1)$ and assume that
$\omega$ is smoothly regularly varying of index $1-\alpha$
at infinity and bounded below by $\omega(t)\ge\eta t^{1-\varepsilon}$
at $0$ for some $\eta>0$ and $\varepsilon\in(0,1)$. Set
%
\begin{equation}
\rho(t) =\bigl(1 + t^2/\omega(t^2)\bigr).\vadjust{\goodbreak}
\end{equation}
Assume that $A$ satisfies~(\ref{MVrho}) and that $\mu$ satisfies
%
\begin{equation}
W(\rho,\mu)=\sup_{s>0}\bigl\{s\mu(\{\rho\circ\delta>s\})\bigr\}<\infty.
\end{equation}
Then we have
\[
\mathcal E_\mu(f,f) \le C_0^2C(\omega)
W(\rho,\mu) \|\psi(A)^{1/2}f\|_2^2,\qquad f\in D_A.
\]
\end{theorem}
\begin{pf} We follow the proof of Theorem
\ref{th-interpol}.
Taking into account that
$\theta(t)\simeq c_\alpha(t/\omega(t))^{1/2}$ and that we have set
$\rho(t)= 1+ t^2/\omega(t^2)$, the proof of Theorem
\ref{th-interpol} shows that we need to estimate
\begin{eqnarray*}
\int_GK(\cdot,s^2)[\rho\circ\delta]^2 \,d\mu&=&
\rho(s)^{1/2}\int_{\{\delta> s \}}
[\rho\circ\delta]^{1/2} \,d\mu\\
&&{}+
\frac{1}{\omega(s^2)^{1/2}}\int_{\{\delta\le s\}}
\delta[\rho\circ\delta]^{1/2} \,d\mu,
\end{eqnarray*}
uniformly over the range $s>1$.
Setting $v(s)= \mu(\delta>s)$, we have
\[
\rho(t)v(t)\le W(\rho,\mu)
\]
and
\begin{eqnarray*}
\rho(s)^{1/2}\int_{\{\delta> s \}}
[\rho\circ\delta]^{1/2} \,d\mu&=&
\rho(s)^{1/2} \int_s^\infty\rho^{1/2}(t) \,d[-v(t)]\\
&\le& \frac{\rho(s)^{1/2}}{2}\int_s^\infty\rho'(t)\rho(t)^{-1/2} v(t)\,dt
+\rho(s)v(s)\\
&\le& W(\rho,\mu)\biggl(1+ \frac{\rho(s)^{1/2} }{2}
\int_s^\infty\rho'(t)\rho(t)^{-3/2}\,dt\biggr)\\
&\le& 2W(\rho,\mu).
\end{eqnarray*}
Further, using the fact that $\omega$ is regularly varying with
positive index $1-\alpha$, we have $t\rho'(t)\sim2\alpha\rho(t)$ and
\begin{eqnarray*}
\frac{1}{\omega(s^2)^{1/2}}\int_{\{\tau\le s\}}
\delta[\rho\circ\delta]^{1/2} \,d\mu&\le&
\frac{1}{\omega(s^2)^{1/2}}\int_0^s
t \rho(t)^{1/2} d[-v(t)] \\
&\le& \frac{1}{\omega(s^2)^{1/2}}\int_0^s
\bigl( \rho(t)^{1/2}+t \rho'(t)\rho(t)^{-1/2} \bigr)v(t)\,dt \\
&\le&
\frac{C(\omega)W(\rho,\mu)}{\omega(s^2)^{1/2}}\int_0^s\rho(t)^{-1/2}\,dt
\\
&\le&
\frac{C(\omega)W(\rho,\mu)}{\omega(s^2)^{1/2}}\int_0^s\frac{\omega
(t^2)^{1/2}\,dt}{t}\\
&\le&
C'(\omega)W(\rho,\mu).
\end{eqnarray*}
This gives the desired result.
\end{pf}
\begin{rem} The case when $\omega$ is a slowly varying increasing function
corresponds to moment conditions that are close to a finite second moment.
In this case, the use of Theorem~\ref{th-interpol} is limited by the
fact that
it involves the possibly infinite quantity
\[
\zeta(t)=t^{1/2}\int_t^\infty\frac{ds}{s\omega(s)^{1/2}}.
\]
We can improve the result by using a slightly different proof.
Namely, using the same notation as in the proof of Theorem~\ref{th-interpol},
we write
\[
g(h)\le2C_0 \int_0^\infty\mathbf K(h,s)([s\omega(s)]^{1/2}
\|AH_{s/2} f\|_2) \,\frac{ds}{s},
\]
where $\mathbf K$ is the kernel on $G\times(0,\infty)$ given by
\[
\mathbf K(h,s)= \frac{s^{1/2}}{\omega(s)^{1/2}}
\bigl(\mathbf1_{(0,\tau(h))}(s)+ \delta(h)s^{-1/2}\mathbf1_{[\tau(h),
\infty)}(s)\bigr).
\]
Next, we use the Hilbert--Schmidt norm
$\int_G\int_0^\infty|\mathbf K(h,s)|^2\,\frac{ds}{s}\,d\mu(h)$
to estimate the norm of $\mathbf K\dvtx L^2((0,\infty),\frac{ds}{s})\ra
L^2(G,d\mu)$.
We have
\begin{eqnarray*}
\int_G\int_0^\infty|\mathbf K(h,s)|^2\,\frac{ds}{s}\,d\mu(h)&=&
\int_G \biggl(\int_0^{\tau(h)}\,\frac{ds}{\omega(s)}
+\delta(h)^2\int_{\tau(h)}^\infty
\,\frac{ds}{s\omega(s)}\biggr) \,d\mu(h)\\
&\le& \int_G \bigl(\widetilde{\xi}^2(\tau(h)) +
\widetilde{\zeta}^2(\tau(h))\bigr)
\,d\mu(h),
\end{eqnarray*}
where
\[
\widetilde{\xi}(t)= \biggl(\int_0^{t}\frac{ds}{\omega(s)}\biggr)^{1/2}
\quad\mbox{and}\quad \widetilde{\zeta}(t)= \biggl(t\int_t^\infty\frac{ds}{s\omega(s)}\biggr)^{1/2}.
\]
This implies\vspace*{1pt} that the conclusion of Theorem~\ref{th-interpol} holds
under the hypothesis that $\rho(t)\ge \widetilde{\zeta}^2(t^2)+
\widetilde{\zeta}^2(t^2)$.

For instance, consider the case when $\omega(t) =[\log(e+t)]^\alpha$,
$\alpha>0$. In this case, $\psi(t)\sim t [\log(e+1/t)]^\alpha$.
On the one hand, we have
$\zeta(t)= \infty$ if $\alpha\le2$ and
$\zeta(t)\simeq t^{1/2}[\log(e+t)]^{1-\alpha/2}$ if $\alpha>2$.
This means that
Theorem~\ref{th-interpol} requires $\alpha>2$ and $
\rho(t) \ge C (1+t)^2 [\log(e+t)]^{1-\alpha}$.

On the other hand, we have
$\widetilde{\zeta}(t)\simeq t^{1/2} [\log(e+t)]^{(1-\alpha)/2}$ if
$\alpha>1$.
This means that the variation explained above requires only $\alpha>1$,
with the same $\rho$, that is,
$\rho(t)\ge C(1+t)^2[\log(e+t)]^{1-\alpha}$.
\end{rem}

\subsection{Two fundamental examples}\label{sec-ex}

\subsubsection*{First example}
Let $G$ be a unimodular Lie group, and let
$(A,D_A)$ be the (unique) self-adjoint extension of a H\"ormander sum
of squares
\[
A=\sum_1^kX_i^2\qquad \mbox{acting on } \mathcal C^\infty_c(G),\vadjust{\goodbreak}
\]
where $\{X_i, i=1,\ldots, k\}$
is a fixed set of left-invariant vector fields which generates the Lie
algebra of $G$
(H\"ormander condition). Then, it is known that $H_tf= f*\mu_t$ where
$(\mu_t)_{t>0}$ is a convolution semigroup of probability measures, and
each $\mu_t$ admits a smooth positive density $x\mapsto h_t(x)$
with respect to the Haar measure $\lambda$; see, for example,
\cite{VSCC}, Chapter 3.
Further, as $t$ tends to infinity,
we have
\begin{eqnarray*}
h_t(e) &\simeq&\Phi_G(t)\\
&\simeq&\cases{
e^{-t}, &\quad if $G$ is not amenable,\vspace*{2pt}\cr
e^{-t^{1/3}}, &\quad if $G$ is amenable with exponential volume
growth,\vspace*{1pt}\cr
t^{-D/2}, &\quad for some integer $D$, otherwise.}
\end{eqnarray*}
For each integer $D$, the last case occurs exactly when $G$ has
polynomial volume growth of degree $D$. The value $h_t(e)$ is the
maximal value of
the function $h_t$ on $G$, and, furthermore, it equals
the norm of the linear operator $H_t\dvtx L^1(G)\ra L^\infty(G)$ as well as the
square of the norm of $H_{t/2}\dvtx L^2(G)\ra L^\infty(G)$. In this case,
we set
\[
\delta(x)=\sup\Biggl\{f(x)-f(e)\dvtx f\in\mathcal C^\infty_c(G),
\sum_1^k|X_if|^2\le1\Biggr\}.
\]
This distance is the sub-Riemannian distance naturally associated with
the set
of left-invariant vector fields $\{X_1,\ldots,X_k\}$, and $\delta(x)$ is finite
for all $x\in G$ because we assume that the $X_i$'s generate the Lie algebra
(this is a special case of one of the fundamental theorem of sub-Riemannian
geometry, often referred to as Chow's theorem); see~\cite{Mont} for a detailed
discussion. Further, it is a simple matter (\cite{VSCC}, Lemma VII.1.1)
to see that
$\mathcal E_A(f,f)=\int\sum_1^k|X_if|^2\,d\lambda$ and
\[
\int|f_h-f|^2\,d\lambda\le\delta(h)^2\int\sum_1^k|X_if|^2\,d\lambda,\qquad f\in
\mathcal C^\infty_c(G), h\in G.
\]
This shows that~(\ref{MVrho}) holds true in this case since $
\int\sum_1^k|X_if|^2\,d\lambda=\|A^{1/2}f\|_2^2$.

\subsubsection*{Second example}

Let $G$ be a compactly generated unimodular group, and set
$Af=f-f*\phi_0$ where $\phi_0$ is continuous, symmetric, compactly supported
probability density on $G$ with the property that $\phi_0>0$ on a
compact generating neighborhood of the identity. Then
\[
\mathcal E_A(f,f)=(1/2)\int_G\|f_h-f\|_2^2\phi_0(h)\,d\lambda(h)
\]
and
$H_tf= f*h_t$ where
\[
h_t=e^{-t}\sum_0^\infty\frac{t^n}{n!} \phi_0^{(n)}.
\]
In particular, if $G$ is a finitely generated group with finite
symmetric generating set $S$
containing the identity, we can set $\phi_0=(\#S)^{-1}\mathbf1_S$.
In any case, for $t\ge1$,
\[
h_t(e)\simeq\phi_0^{(2t)}(e)\simeq\Phi_G(t).
\]
As explained in the \hyperref[sec-Intr]{Introduction}, many different behaviors are possible for
the function $\Phi_G$, depending on $G$. Assuming that $U$ is a symmetric
neighborhood of the identity which contains a generating compact set, that
$\inf_{U^3}\{\phi_0\}>0$, and setting
\[
\delta(x)=\inf\{n\dvtx x\in U^n\},
\]
\cite{VSCC}, Proposition VII.3.2, gives that (the discrete case
of this inequality is a bit simpler and more elementary)
\[
\int|f_h-f|^2\,d\lambda\le C(U,\phi_0) \delta(h)^2 \mathcal
E_A(f,f),\qquad
f\in L^2(G), h\in G.
\]
Again, this shows that~(\ref{MVrho}) holds true in this setting.

\section{\texorpdfstring{Applications: Main lower bounds on $\Phi_{G,\rho}$}
{Applications: Main lower bounds on Phi G, rho}}\label{sec-low-appl}

Let $G$ be as in the second example of Section~\ref{sec-ex}.
Keep the notation introduced there. In the applications we have in mind,
we are given a continuous increasing function $\rho\dvtx(0,\infty)\ra
[1,\infty)$
and set $\rho_G=\rho\circ\delta$.
Our main aim is to estimate the functions $\Phi_{G,\rho_G}$ and
$\widetilde{\Phi}_{G,\rho_G}$ introduced in
Definitions~\ref{def1} and~\ref{def2}.
Hence, we consider a (otherwise arbitrary)
symmetric continuous probability density $\phi$ on $G$ with the
property that
$\|\phi\|_\infty\le K$ and
$\int\rho_G \phi \,d\lambda\le K \sup_{\Omega^2}\{\rho\}$ or
$W(\rho_G,\phi)\le K \sup_{\Omega^2}\{\rho\}$.
Here $K>1$ and $\Omega$ are as in Definitions~\ref{def1} and~\ref{def2}.

In order to apply Theorem~\ref{th-interpol},
we have to find an increasing function $\omega$ compatible with $\rho$
in the
sense that the pair $\rho,\omega$ satisfies
the various hypotheses of Theorem~\ref{th-interpol}. The function $\psi$
associated to $\rho$ via $\omega$ is then defined by~(\ref{omegapsi}).

The following examples
are of particular interest:

\begin{itemize}
\item If $\rho(s)=\rho_{2\alpha}(s)=(1+s)^{2\alpha}$, $\alpha\in(0,1]$,
then we can take
\[
\omega(s)=\Gamma(2-\alpha)^{-1} s^{1-\alpha} \quad\mbox{and}\quad
\psi(s)=s^\alpha.
\]
\item If $\rho(s)=(1+s^2)^{\alpha} \ell(1+s^2)^{\alpha} $
with $\alpha\in(0,1)$ and $\ell$ smooth, positive and slowly varying
at infinity, then we can take
\[
\omega(s)= \frac{1+s}{[(1+s)\ell(1+s)]^\alpha} \qquad\mbox{at infinity}
\]
and
\[
\psi(s)\simeq[(1+s)/\ell(1+1/s)]^\alpha\qquad\mbox{at }0.
\]
\item If $\rho(s)= \rho^{\exp}_{c,\alpha}(s)=\exp( c[\log(1+s)]^\alpha)$,
$\alpha\in(0,1)$, $c>0$, then
we can take
\[
\omega(s)= s \exp\bigl( -c_1[\log(1+s)]^\alpha\bigr)\vadjust{\goodbreak}
\]
for some $c_1>0$ (see Remark~\ref{rem-interpol})
and we have [see~(\ref{omegapsi})]
\[
\psi(s)\sim \exp\bigl( -c_1[\log(1+1/s)]^\alpha\bigr) \qquad\mbox{at } 0.
\]
\item If $\rho(s)= \rho^{\log}_\alpha(s)= [\log(e+s)]^\alpha$,
$\alpha>1$, then
we can take
\[
\omega(s)= s[\log(e+s)]^{1-\alpha}
\]
for some $c_1>0$ (see Remark~\ref{rem-interpol}), and this gives
\[
\psi(s)\sim[\log(e+1/s)]^{1-\alpha} \qquad\mbox{at } 0.
\]
These computations indicate that
Theorem~\ref{th-interpol} is too weak to provide results when
$\rho(s)= \rho^{\log}_\alpha(s)= [\log(e+s)]^\alpha$ and $\alpha\le1$.

\item The previous two cases can be generalized as follows. Assume
that
\[
\rho(s)\ge c\log(e+s) \bigl(1+\eta(s^2)\bigr),
\]
where $\eta$ is a positive increasing function such that
$\eta(s) \sim s$ at $0$ and $\eta$ is slowly varying at infinity.
Set
\[
\omega(s) = s/\eta(s) \quad\mbox{and}\quad
\psi(\lambda)= \lambda^2 \int_0^\infty e^{-\lambda s}\omega(s)\,ds.
\]
By~\cite{BGT}, Theorem 1.7.1, we have
$\psi(\lambda)\simeq c/\eta(1/\lambda)$. Further, referring
to the notation used in Theorem~\ref{th-interpol}, we then have
\[
\frac{t\max\{\xi(t^2),\zeta(t^2)\}}{\omega(t^2)^{1/2}}\le C\rho(t).
\]
\end{itemize}

We are now ready to state and prove lower bounds on
the functions $\Phi_{G,\rho_G}$ and $\widetilde{\Phi}_{G,\rho_G}$
of Definitions~\ref{def1} and~\ref{def2} for some
groups $G$ and functions $\rho$. We will use the notation
$\rho_\alpha,\rho^{\exp}_{c,\alpha},\rho^{\log}_\alpha$ recalled above.
If $\rho$ is a real function, and $G$ is a compactly generated group with
world length
$\delta(x)=|x|=\inf\{n\dvtx x\in U^n\}$ for some fixed symmetric relatively compact
generating neighborhood of the identity, we set
$\Phi_{G,\rho}=\Phi_{G,\rho_G}$ where $\rho_G=\rho\circ\delta$.

We state four theorems that cover various cases of particular interest.
The proofs of these results all follow the same outline based on
Theorems~\ref{th-interpol} and~\ref{th-interpol2} together with
Corollary~\ref{cor-comp1} and the results of Section
\ref{sec-funct}, Theorems~\ref{th-RRP} and~\ref{th-RRR}.
The main line of reasoning described in the proof of Theorem
\ref{th-Phi>1} below is also used for the proofs of Theorems~\ref{th-Phi>2},
\ref{th-Phi>3} and~\ref{th-Phi>4}. The results presented in the
\hyperref[app]{Appendix}
play a crucial role in these proofs.
\begin{theorem} \label{th-Phi>1}
Let $G$ be a locally compact, compactly generated unimodular group
such that $\Phi_G(n)\simeq n^{-D/2}$ at infinity, for some integer $D$.
\begin{longlist}[(2)]
\item[(1)] Assume that
$\rho(s)\ge[(1+ s^2)\ell(1+s^2)]^\alpha$ with $\alpha\in(0,1)$
and $\ell$ smooth positive slowly varying at infinity with de Bruijn conjugate
$\ell^{\#}$. Then there exist $c=c_\rho\in(0,\infty)$ and an integer
$N=N_\rho$ such that
\[
\forall n>N \qquad \Phi_{G,\rho}(n)\ge\widetilde{\Phi}_{G,\rho}(n)\ge
c [n^{1/\alpha}\ell^{\#}(n^{1/\alpha})]^{-D/2}.
\]
\item[(2)] Assume that $\rho(s)\ge \log(e+s) \ell(1+s^2)$ and $\ell$
smooth positive increasing and slowly varying at infinity and such that
$\log\ell^{-1}(t) \simeq t^\gamma\omega(t)^{1+\gamma}$ at infinity,
with $\gamma\ge0$ and $\omega$ slowly varying with de Bruijn conjugate
$\omega^{\#}$. Then there exist $C=C_\rho\in(0,\infty)$ and an integer
$N=N_\rho$
such that
\[
\forall n>N\qquad
\log\Phi_{G,\rho}(n)\ge
-C [n^\gamma/\omega^{\#}(n)]^{1/(1+\gamma)}.
\]
\end{longlist}
\end{theorem}
\begin{pf}
We will use the following notation which is consistent with the
notation used in the \hyperref[app]{Appendix}.
Let $\psi\dvtx[0,2]\ra[0,2]$ be a continuous increasing
function with continuous derivative such that
with $\psi(0)=0$, $\psi(1)=1$ and $\psi(2)<2$.
Fix a symmetric probability density $\phi_0\in L^2(G)$ and assume that
its support is a compact generating neighborhood of the identity element
(we assume that $G$ is compactly generated).
This implies that $\phi^{(2n)}_0(e)\simeq
\Phi_G(n)$; see~\cite{Gre,PSCstab,VSCC}.
We set $T=R_{\phi_0}$ and [see~(\ref{Tpsi})]
\[
T_\psi= I- \psi(I-T).
\]

Let $\mathcal E_0$
denote Dirichlet form
$\mathcal E_0(f,f)=\langle(I-T)f,f\rangle$ associated with $\phi_0$.
By Section~\ref{sec-ex} and Theorems~\ref{th-interpol} and~\ref{th-interpol2},
if $d\mu=\phi \,d\lambda$ is a symmetric probability with continuous density
satisfying $\mu(\rho\circ\delta)<\infty$, then
%
\begin{equation}\label{cruxcomp}
\mathcal E_\mu(f,f)\le C\mathcal\| \psi(I-T)^{1/2}f\|_2^2.
\end{equation}
Here $\psi$ is chosen such that the condition of
Theorems~\ref{th-interpol} and~\ref{th-interpol2} relating $\psi$ to $\rho$
(via $\omega$) are satisfied. See the explicit examples discussed at the
beginning of this section.

By Corollary~\ref{cor-comp1},~(\ref{cruxcomp}) implies
[$\tau$ is the natural semifinite trace on the
von Neumann $V(G)$; see the \hyperref[app]{Appendix}]
\[
\phi^{(2n)}(e) \ge C\bigl(e^{-cn}+ \tau\bigl(T_\psi^{2[cn]}\bigr)\bigr).
\]
Now, depending on the behavior of $\psi$ near $0$,
the trace $\tau(T_\psi^{2n})$ can be
estimated using the results of Section
\ref{sec-funct}, Theorems~\ref{th-RRP} and~\ref{th-RRR}; see also
Example~\ref{exa-gamma-omega}. This gives the announced lower bounds
on $\phi^{(2n)}(e)$.
\end{pf}
\begin{exa} The second statement in Theorem~\ref{th-Phi>1}
can be illustrated by the following two examples:
\begin{longlist}[(2)]
\item[(1)] $\log\Phi_{G,\rho^{\exp}_{c,\alpha}}(n) \ge-
C_{D,c,\alpha} (\log n)^{1/\alpha}$,
$\alpha\in(0,1)$, $c>0$.
\item[(2)]
$\log\Phi_{G,\rho^{\log}_\alpha}(n) \ge- C_{D,\alpha} n^{1/\alpha}$,
$\alpha>1$.
\end{longlist}
\end{exa}
\begin{theorem} \label{th-Phi>2}
Let $G$ be a locally compact, compactly generated unimodular group
such that $\log\Phi_G(n) \ge-C n^{\gamma}$ at infinity, for some
$\gamma\in(0,1)$ and $C\in(0,\infty)$. Fix $\alpha\in(0,1)$ and
$\rho(s)\simeq[(1+s^2)\ell(1+s^2)]^\alpha$ with
$\ell$ smooth positive slowly varying at infinity. Then,
there exist $C_\rho\in(0,\infty)$ and an integer $N_\rho$
such that
\[
\forall n>N_\rho\qquad
\log\Phi_{G,\rho}(n)\ge\log\widetilde{\Phi}_{G,\rho}(n)\ge -C_\rho
n^{\gamma_\alpha}/\ell_*^{\#}\bigl(n^{(1-\gamma)\gamma_\alpha/\gamma}\bigr)
^{\alpha},\vadjust{\goodbreak}
\]
where $\ell^{\#}_*$ is the de Bruijn conjugate of
\[
\ell_*(s)= \ell^{\#}(s)^{\gamma_\alpha},\qquad
\gamma_\alpha=\frac{\gamma}{\gamma+\alpha(1-\gamma)}.
\]
\end{theorem}

This theorem with $\ell\equiv1$ implies Theorem~\ref{th-I2}.
\begin{pf*}{Proof of Theorem~\ref{th-Phi>2}}
The given $\rho$ calls for using $\psi(s)= [s/\ell(s)]^\alpha$
in Theorem~\ref{th-interpol}. Note that $\psi^{-1}(t)
\simeq t^{1/\alpha}/\ell^{\#}(1/t^{1/\alpha})$.

Using Theorem~\ref{th-RRR} and
the same notation and line of reasoning as in the proof of Theorem
\ref{th-Phi>1}, we obtain that
if $d\mu=\phi \,d\lambda$ is a symmetric probability with continuous density
satisfying $\mu(\rho\circ\delta)<\infty$, then
\[
\log\phi^{(2n)}(e) \ge-C_1n/\pi_{\psi}(n)
\]
with
\[
C\pi_\psi^{-1}(Ct)\ge t \psi^{-1}(1/t)^{-\gamma/(1-\gamma)} \ge c
t^{(\alpha(1-\gamma)+\gamma)/\alpha(1-\gamma)}
\ell^{\#}(t^{1/\alpha})^{\gamma/(1-\gamma)}.
\]
This can be written as (for a different constant $C$)
\[
C\pi_\psi^{-1}(Ct)\ge t^{(\alpha(1-\gamma)+\gamma)/\alpha(1-\gamma)}
\ell_*(t^{1/\alpha})^{(\alpha(1-\gamma)+\gamma)/(1-\gamma)}
\]
with $\ell_*(s)= \ell^{\#}(s)^{\gamma/(\alpha(1-\gamma)+\gamma)}$.
This gives
\[
c\pi_\psi(ct)\le t^{\alpha(1-\gamma)/(\alpha(1-\gamma)+\gamma)}
\ell_*^{\#}\bigl(t^{(1-\gamma)/(\alpha(1-\gamma)+\gamma)}\bigr)^\alpha
\]
and
\[
\log\Phi_{G,\rho}(n) \ge-Cn^{\gamma_\alpha}/\ell_*^{\#}\bigl(n^{(1-\gamma
)/(\alpha(1-\gamma)+\gamma)}\bigr)^\alpha
\]
with
$ \gamma_\alpha=\gamma/(\alpha(1-\gamma)+\gamma)$, as desired.
\end{pf*}
\begin{exa} Assume that $\ell$
satisfies $\ell(t^a)\simeq\ell(t)$ for all $a>0$. Then $\ell^{\#
}=1/\ell$ and
$(1/\ell)^{\#}\simeq\ell$. Hence
$\ell_*\simeq(1/\ell)^{\gamma/(\alpha(1-\gamma)+\gamma)}$ and
$\ell_*^{\#}\simeq\ell ^{\gamma/(\alpha(1-\gamma)+\gamma)}$. Hence
we get
\[
-\log\Phi_{G,\rho}(n) \le C [n/\ell^\alpha(n)]^{\gamma_\alpha}.
\]
This is consistent with Example~\ref{exa-gamma}.
\end{exa}
\begin{theorem} \label{th-Phi>3}
Let $G$ be a locally compact, compactly generated unimodular group
such that $\log\Phi_G(n) \ge-C n^{\gamma}/\ell(n)$ at infinity, for some
$\gamma\in(0,1]$, $C\in(0,\infty)$ and slowly varying function $\ell$
satisfying $\ell(t^a)\simeq\ell(t)$ for all $a>0$.
Assume that $\alpha\in(0,1)$ and
$\rho(s)= (1+s)^{2\alpha}$. Then, we have
\[
\log\Phi_{G,\rho}(n)\ge\log\widetilde{\Phi}_{G,\rho}(n)\ge
-C_\rho[n/\ell(n)^{\alpha/\gamma}]^{\gamma_\alpha},\qquad
\gamma_\alpha=\frac
{\gamma}{\gamma+\alpha(1-\gamma)} .
\]
\end{theorem}
\begin{theorem}\label{th-Phi>4}
Let $G$ be a locally compact, compactly generated unimodular group
such that $\log\Phi_G(n) \ge-C n/\pi(n)$ with $\pi$ continuous increasing
and satisfying $\pi(t)\le t^{1-\varepsilon}$ for $t\ge1$.
Assume that $\rho(s)\ge c \log(e+s) \ell(s^2)$ with $c>0$\vadjust{\goodbreak} and $\ell$
smooth positive increasing and slowly varying at infinity.
Then there exist $c_1,C_1\in(0,\infty)$ such that, for all $n$ large enough,
\[
\log\Phi_{G,\rho}(n) \ge - C_1 n/\ell(\pi(c_1 n^\varepsilon)).
\]
\end{theorem}
\begin{exa} If $\Phi_G(n)\ge c\exp(-Cn^\gamma)$ with $\gamma\in(0,1)$,
this yields:
\begin{longlist}[(2)]
\item[(1)] $\log\Phi_{G,\rho^{\exp}_{c,\alpha}}(n) \ge -
C_{\gamma,c,\alpha} n\exp(-c_{\gamma,c,\alpha}[\log n]^\alpha) $,
$\alpha\in(0,1)$, $c>0$.
\item[(2)] $\log\Phi_{G,\rho^{\log}_\alpha}(n) \ge -
C_{\gamma,\alpha}n[\log n]^{-(\alpha-1)}$,
$\alpha>1$.
\end{longlist}

If, instead, $\Phi_G(n)\ge c\exp(-Cn/\ell(n))$ with $\ell$
increasing slowly varying and satisfying $\ell(t^a)\simeq\ell(t)$ for
all $a>$,
we obtain:
\begin{longlist}[(2)]
\item[(1)] $\log\Phi_{G,\rho^{\exp}_{c,\alpha}}(n) \ge -
C_{\gamma,c,\alpha} n\exp(-c_{\gamma,c,\alpha}[\log\ell(n)]^\alpha) $,
$\alpha\in(0,1)$, $c>0$.
\item[(2)] $\log\Phi_{G,\rho^{\log}_\alpha}(n) \ge -
C_{\gamma,\alpha}n[\log\ell(n)]^{-(\alpha-1)}$,
$\alpha>1$.
\end{longlist}
\end{exa}
\begin{pf*}{Proof of Theorems~\ref{th-Phi>3} and~\ref{th-Phi>4}}
In each case, we use either Theorem~\ref{th-interpol} or
Theorem~\ref{th-interpol2}
together with either Theorems~\ref{th-RRP} or~\ref{th-RRR},
and Corollary~\ref{cor-comp1}.
\end{pf*}

\section{\texorpdfstring{Upper bounds on $\Phi_{G,\varrho}$}
{Upper bounds on Phi G, rho}}

The aim of this section is to obtain upper bounds on the function $\Phi
_{G,\varrho}$
(and its variant $\widetilde{\Phi}_{G,\varrho}$) under various conditions
on the group $G$ and the function $\varrho$. To obtain such upper bounds,
we only need to exhibit an example of a symmetric probability density
$\phi$ such that
$\int\varrho\phi \,d\lambda<\infty$ (or $\sup_{s>0}\{s\int_{\{\varrho
>s\}}\phi \,d\lambda\}$,
in the case of $\widetilde{\Phi}_{G,\varrho}$)
and for which we can obtain an upper bound on
$n\mapsto\phi^{(2n)}(e)$. Of course, to obtain good upper bounds, we
need to
identify probability densities with the desired moment condition and
for which $n\mapsto\phi^{(2n)}(e)$ presents an almost optimal decay.
This question---which densities produce the optimal
decay?---is quite interesting in its own right. For instance, when $G$ is finitely
generated with finite symmetric generating set $S$ and $\varrho$ is of
the form
$\varrho=\rho_G=\rho(|\cdot|)$, and $|\cdot|=|\cdot|_S$ is the
word-length based on the
generating set $S$, should we expect to find a probability density with
nearly optimal decay among ``radial densities'' of the form $\phi(x)= f(|x|)$?

\subsection{\texorpdfstring{$\Phi_G$-based upper bounds: subordination}
{Phi G-based upper bounds: subordination}} \label{sec-Sub}

The lower bounds on $\Phi_{G,\varrho}$ (and $\widetilde{\Phi}_{G,\varrho
}$) obtained
in Section~\ref{sec-low-appl} for certain $\varrho=\rho\circ\delta$
are all based on lower bounds on the function $\Phi_G$. It is thus natural
to seek upper bounds of the same nature.
These applications of Theorem~\ref{th-interpol}
start with a symmetric compactly supported continuous density
$\phi$ (with generating support) and involve comparison with the behavior
of certain operators $T_\psi$ of the form $T_\psi=I -\psi(I-R_\phi)$
where the function $\psi$ is chosen appropriately, depending on $\rho$.

It would be very nice to
identify a class of functions $\psi$ so that $T_\psi= R_{\phi_\psi}$
where $\phi_\psi$ is, itself, a symmetric probability density. As already
noted after~(\ref{Tpsi}), this is certainly the case when $\psi$ is a
Bernstein function satisfying\vadjust{\goodbreak} $\psi(0)=0$, $\psi(1)=1$; see, for
example,~\cite{NJ,SSV} for an access to the literature on
Bernstein functions.
A Bernstein function is a smooth positive function
$\psi\dvtx(0,\infty)\ra(0,\infty)$ such that $(-1)^k\,\frac{d^k\psi}{dt^k}
\le0$
and two\vspace*{-2pt} good and important examples of Bernstein functions are
$\psi_\alpha\dvtx s\mapsto s^\alpha$,
$\alpha\in(0,1]$ and
$\psi^{\log}_\alpha\dvtx s\mapsto[\log_2(1 +s^{-1/\alpha})]^{-\alpha}$.
Further, for any smooth positive increasing regularly varying function
$\psi_1$ of index $\alpha$ in $[0,1)$ at $0$ such that
$x\mapsto x\psi'_1(x)$ is also regularly varying of index $\alpha$,
there exists a Bernstein function $\psi$ such that
$\psi\sim\psi_1$; see~\cite{BSCsubord}, Theorem 2.5.

In order to obtain upper bounds on the functions
$\Phi_{G,\rho_G}$ and $\widetilde{\Phi}_{G,\rho_G}$,
it suffices to find a Bernstein\vspace*{-1pt} function $\psi$ such that the probability
density $\phi_\psi$ satisfies the required moment condition and to estimate
$\phi_\psi^{(2n)}(e)$. The companion paper~\cite{BSCsubord} develops this
idea, and we will simply quote the relevant results.

We start with results concerning groups with polynomial volume growth
$V(n)\simeq n^D$. By~\cite{HSC}, these groups satisfy $\Phi_G(n)\simeq
n^{-D/2}$.
In fact, thanks to~\cite{HSC} and deep results of Guivarc'h, Gromov and Losert,
groups of polynomial volume growth are exactly those groups that satisfy
$\Phi_G(n)\simeq n^{-D/2}$ for some integer $D$.
An alternative and
self-contained proof of
the theorems discussed below is given in the next section.
\begin{theorem}[(\cite{BSCsubord})]\label{th-Phi-Pol1}
Assume\vspace*{1pt} that $G$ is a compactly generated locally compact
group with polynomial volume growth $V(n)\simeq n^D$.
\begin{longlist}[(2)]
\item[(1)] Assume that $\rho(s)\simeq g(1+s^2)$ where $g(s)= [s \ell
(s)]^\alpha$
where $\alpha\in(0,1)$, and
$\ell$ is a positive slowly varying function at infinity with
de Bruijn conjugate $\ell^{\#}$.
Then there exist $C\in(0,\infty)$ and an integer $N$ such that
\[
\forall n>N\qquad
\widetilde{\Phi}_{G,\rho_G}(n)\le
C[n^{1/\alpha}\ell^{\#}(n^{1/\alpha})]^{-D/2}.
\]
Further, for any slowly varying function $\ell_1$ with de Bruijn
conjugate $\ell_1^{\#}$
such that $\sum_1^\infty\frac{\ell(n)^\alpha}{n \ell_1(n)^\alpha
}<\infty$,
there exist $C(\ell_1)\in(0,\infty)$ and an integer $N(\ell_1)$ such that
\[
\forall n>N(\ell_1)\qquad
\Phi_{G,\rho_G}(n)\le
C(\ell_1)[n^{1/\alpha}\ell_1^{\#}(n^{1/\alpha})]^{-D/2}.
\]
\item[(2)] Assume that\vspace*{1pt} $\rho(t)\simeq g(1+s^2)$ where $g (s)= \widehat
{\ell}(s)$
and $\widehat{\ell}(t)= 1/\int_t^\infty\frac{du}{u \ell(u)}$ where $\ell$
is a
positive increasing slowly varying function at infinity. Assume further that
$\log\widehat{\ell}^{-1}(t)\simeq t^\gamma\omega(t)^{1+\gamma}$ at infinity,
with $\gamma\ge0$ and $\omega$ slowly varying with de Bruijn conjugate
$\omega^{\#}$. Then there exist $c,C\in(0,\infty)$ and an integer $N$
such that
\[
\forall n>N\qquad
\widetilde{\Phi}_{G,\rho_G}(n)\le
C\exp\bigl(-c [n^\gamma/\omega^{\#}(n)]^{1/(1+\gamma)}\bigr).
\]
Further, for any slowly varying function $\ell_1$ such that
\[
\sum_1^\infty\frac{\widehat{\ell}(n)}{n \ell_1(n)}<\infty
\quad\mbox{and}\quad
\log\widehat{\ell}_1^{-1}(t)\simeq t^{\gamma_1}\omega_1(t)^{1+\gamma_1}\vadjust{\goodbreak}
\]
with $\gamma_1\ge0$ and $\omega_1$ slowly varying at infinity,
there exist $c=c(\ell_1), C=C(\ell_1)\in(0,\infty)$ and an integer
$N=N(\ell_1)$ such that
\[
\forall n>N\qquad
\Phi_{G,\rho_G}(n)\le
C\exp\bigl(-c [n^{\gamma_1}/\omega_1^{\#}(n)]^{1/(1+\gamma_1)}\bigr).
\]
\end{longlist}
\end{theorem}

Putting together the results of Theorems~\ref{th-Phi>1} and~\ref{th-Phi-Pol1},
we obtain the following results which imply Theorem~\ref{th-Ipol}.
\begin{theorem}Assume that $G$ is a compactly generated locally compact
group with polynomial volume growth $V(n)\simeq n^D$.
\begin{longlist}[(3)]
\item[(1)] Assume that $\rho(s)\simeq g(1+s^2)$ where $g(s)= [s \ell
(s)]^\alpha$
where $\alpha\in(0,1)$ and
$\ell$ is a positive slowly varying function at infinity with
de Bruijn conjugate $\ell^{\#}$. Then
\[
\widetilde{\Phi}_{G,\rho_G}(n)\simeq[n^{1/\alpha}\ell^{\#}(n^{1/\alpha
})]^{-D/2}.
\]


\item[(2)] For any $\alpha\in(0,1)$ and $c>0$, there are constants
$c_1,c_2,C_1,C_2$
(depending on $G$, $\alpha$ and $c$) such that
\[
\forall n\qquad c_1 \exp( -C_1 [\log n]^{1/\alpha}) \le
\Phi_{G,\rho^{\exp}_{c,\alpha}}(n) \le
C_2\exp(-c_2[ \log n]^{1/\alpha} ).
\]

\item[(3)] For any $\beta>\alpha>1$, there are constants $c_1,c_2,C_1,C_2$
(depending on $G$, $\alpha$ and $\beta$) such that
\[
\forall n\qquad c_1 \exp(-C_1 n ^{1/\alpha}) \le
\Phi_{G,\rho^{\log}_\alpha}(n) \le C_2\exp\bigl(-c_2 n^{1/(\beta+1)}
\bigr).
\]
\end{longlist}
\end{theorem}

Our next result concerns groups with volume growth faster than
polynomial and moment of the type $\rho_\alpha(s)=(1+s)^{\alpha}$.
No classifications of either volume growth or the behavior of $\Phi_G$
are known
for such groups. The upper bounds in the following theorem cannot
be obtained by the methods
developed in the next section.
This theorem follows immediately from Theorems~\ref{th-Phi>2} and~\ref{th-Phi>3}
and~\cite{BSCsubord}, Theorem 5.3.
\begin{theorem}\label{th-superpol}
Assume that $G$ is a compactly generated locally compact
group and that there exist $0\le\bar{\gamma}\le\gamma\le1$ and positive
slowly varying functions $\eta,\bar{\eta}$, both satisfying $\eta
(t^a)\simeq
\eta(t)$ for all $a>0$, such that, for $n$ large enough,
\[
- n^{\gamma}/\eta(n) \le \log\Phi_G(n)\le
-n^{\bar{\gamma}}/\bar{\eta}(n).
\]
For any $s>0$, set
$\gamma_s= \gamma/[s(1-\gamma)+\gamma]$,
$\bar{\gamma}_s= \bar{\gamma}/[s(1-\bar{\gamma})+\bar{\gamma}]$.

Assume further
that there exists a symmetric
continuous probability density $\phi$ with compact support,
positive on a generating
compact set and such that
%
\begin{equation}\label{escape-theta}
\forall n,s\qquad \int\mathbf1_{\{|x|\ge n^\theta s\}}
\phi^{(n)}\,d\lambda\le C\exp(-cs^q)
\end{equation}
for some $C,\theta,q>0$.

\begin{longlist}[(2)]
\item[(1)] For any $\alpha\in(0,\min\{2,1/\theta\})$,
there exist $c_1,C_1\in(0,\infty)$ such that, for all $n$ large enough,
\[
-C_1 [n/\eta(n)^{\alpha/2\gamma}]^{\gamma_{\alpha/2}}
\le\log\widetilde{\Phi}_{G,\rho_\alpha}(n) \le - c_1
[n/\bar{\eta}(n)^{\alpha\theta/\bar{\gamma}}]^{\bar{\gamma}_{\alpha
\theta}} .
\]

\item[(2)] For any $\alpha\in(0,\min\{2,1/\theta\})$ and $\varepsilon>0$,
there exist $c_1,C_1\in(0,\infty)$ such that, for all $n$ large enough,
\[
-C_1 [n/\eta(n)^{\alpha/2\gamma}]^{\gamma_{\alpha/2}}
\le\log\Phi_{G,\rho_\alpha}(n) \le - c_\varepsilon
\bigl[n/[\bar{\eta}(n)(\log n)^{1+\varepsilon}]^{\alpha\theta/\bar{\gamma
}}\bigr]^{\bar{\gamma}_{\alpha\theta}} .
\]
\end{longlist}
\end{theorem}
\begin{exa} Assume that $G= F\wr H$ where $F$ is a nontrivial
finite group, and $H$ is polycylic with exponential volume growth. Then
$\Phi_{G}(n)\simeq\exp(-n/(\log n)^2)$.
Condition~(\ref{escape-theta}) is trivially verified with $\theta=1$.
For $\alpha\in(0,2)$, Theorem~\ref{th-superpol}(1) yields
\[
- C_1 n/ [\log n]^{\alpha}\le \Phi_{G,\rho_\alpha}(n)
\le-c_1 n/[\log n]^{2\alpha}
\]
for all $n$ large enough. We conjecture that the lower bound is correct.
\end{exa}

We now state two corollaries of Theorem~\ref{th-superpol}.
The first corollary gives a result that
is widely applicable
whereas the second corollary requires a precise understanding
of the most basic random walks on the group $G$. In particular,
the hypothesis~(\ref{escape-root}) made in Corollary~\ref{cor-superpol2}
requires a classical $\sqrt{n}$ rate of escape for simple random walk
on $G$.
\begin{cor} \label{cor-superpol1}
Assume that $G$ is a compactly generated locally compact
group and that there exist $0\le\bar{\gamma}\le\gamma\le1$ such that,
for $n$ large enough,
\[
- n^{\gamma} \le \log\Phi_G(n)\le
-n^{\bar{\gamma}}.
\]
For any $s>0$, set
$\gamma_s= \gamma/[s(1-\gamma)+\gamma]$,
$\bar{\gamma}_s= \bar{\gamma}/[s(1-\bar{\gamma})+\bar{\gamma}]$.

\begin{longlist}[(2)]
\item[(1)] For any $\alpha\in(0,1)$, there exist $c_1,C_1\in(0,\infty)$
such that, for all $n$ large enough,
\[
-C_1 n^{\gamma_{\alpha/2}}
\le\log\widetilde{\Phi}_{G,\rho_\alpha}(n) \le - c_1
n^{\bar{\gamma}_{\alpha}} .
\]
\item[(2)] For any $\alpha\in(0,1)$ and $\varepsilon>0$,
there exist $c_\varepsilon,C_1\in(0,\infty)$ such that, for all $n$ large enough,
\[
-C_1 n^{\gamma_{\alpha/2}}
\le\log\Phi_{G,\rho_\alpha}(n) \le - c_\varepsilon
n^{\bar{\gamma}_\alpha}/(\log n)^{(1+\varepsilon){\bar{\gamma}_\alpha
}\alpha/\bar{\gamma}} .
\]
\end{longlist}
\end{cor}
\begin{cor}\label{cor-superpol2}
Assume that $G$ is a compactly generated locally compact
group and that there exist $\gamma\in(0,1)$ and a positive
slowly varying function $\eta$ satisfying $\eta(t^a)\simeq
\eta(t)$ for all $a>0$, such that
\[
\log\Phi_G(n)\simeq
-n^{\gamma}/\eta(n).
\]
For any $s>0$, set
$\gamma_s= \gamma/[s(1-\gamma)+\gamma]$.
Assume further
that there exists a symmetric
continuous probability density $\phi$ with compact support,\vadjust{\goodbreak}
positive on a generating
compact set, and such that
%
\begin{equation}\label{escape-root}
\forall n,s\qquad \int\mathbf1_{\{|x|\ge n^{1/2} s\}}
\phi^{(n)}\,d\lambda\le C\exp(-cs^q)
\end{equation}
for some $C,q>0$.
\begin{longlist}[(2)]
\item[(1)] For any $\alpha\in(0,2)$, there exist $c_1,C_1\in(0,\infty)$
such that, for all $n$ large enough,
\[
-C_1 [n/\eta(n)^{\alpha/2\gamma}]^{\gamma_{\alpha/2}}
\le\log\widetilde{\Phi}_{G,\rho_\alpha}(n) \le - c_1
[n/\eta(n)^{\alpha/2\gamma}]^{\gamma_{\alpha/2}} .
\]

\item[(2)] For any $\alpha\in(0,2)$ and $\varepsilon>0$,
there exist $c_\varepsilon,C_1\in(0,\infty)$ such that, for all $n$ large enough,
\[
-C_1 [n/\eta(n)^{\alpha/2\gamma}]^{\gamma_{\alpha/2}}
\le\log\Phi_{G,\rho_\alpha}(n) \le - c_\varepsilon
\bigl[n/[\eta(n)(\log n)^{1+\varepsilon}]^{\alpha/2\gamma}\bigr]^{\gamma_{\alpha/2}} .
\]
\end{longlist}
\end{cor}
\begin{exa} Let the group $G$ be either the group $\mathrm{Sol}=\mathbb
Z\ltimes_A \mathbb Z^2$ where $A=\bigl({2\atop1}\enskip{1\atop1}\bigr)$,
or the wreath product $F\wr\mathbb Z$ where $F$ is any finite group.
By~\cite{Revesc}, these groups satisfy~(\ref{escape-root}). Further,
these groups have exponential volume growth and satisfy $\Phi
_G(n)\simeq\exp(-n^{1/3})$; see, for example,~\cite{VSCC} and the
references therein. Hence Corollary
\ref{cor-superpol2} applies. In particular,
for any $\alpha\in(0,2)$, we have
\[
\widetilde{\Phi}_{G,\rho_\alpha}(n) \simeq
\exp\bigl( - n^{1/(1+\alpha)}\bigr).
\]
Using a different argument,
we shall see in the next section that this result also holds
for all polycyclic groups.
\end{exa}

The final two results of this section concern groups with
super-polynomial volume growth and slowly varying moment condition.
\begin{theorem}\label{th-superpol-low1}
Assume that $G$ is a compactly generated locally compact
group and that there exist $0< \bar{\gamma}\le\gamma < 1$
and $c,C\in(0,\infty)$ such that, for $n$ large enough,
\[
- Cn^{\gamma} \le \log\Phi_G(n)\le
-c n^{\bar{\gamma}}.
\]
Let $\rho(t)\simeq\log(e+t)\ell(t)$ where $\ell$ is a continuous increasing
slowly varying function at infinity.
Let $\rho_1$ be a slowly varying function such that
$\sum_1^\infty\frac{\rho(n)}{n\rho_1(n)}<\infty$,
set $\widehat{\rho}_1(t)=1/\int_t^\infty\frac{ds}{s\rho_1(s)}$
and fix $\varepsilon\in(0,1)$.
Then there are $C_1(\varepsilon),c_1(\rho_1)\in(0,\infty)$ such that,
for all $n$ large enough,
\[
-C(\varepsilon) n/\ell(n^{\varepsilon\gamma/2})
\le \log\Phi_{G,\rho}(n) \le- c(\rho_1) n/\widehat{\rho}_1(n^{\bar
{\gamma}}).
\]
\end{theorem}

Theorem~\ref{th-superpol-low1} applies to a very large collection of groups.
For instance, it applies to all polycyclic groups with exponential volume
growth since such groups have $\Phi_{G}(n)\simeq\exp(-n^{1/3})$.
It also applies to groups with volume growth satisfying
$ c n^{a} \le\log V(n)\le C n^{b}$ with $0<a\le b<1$ since these volume
estimates imply
$ -C_1n^b\le \log\Phi_{G}(n)\le-c_1 n^{a/(a+2)}$.\vadjust{\goodbreak}

The following two examples provide a proof of the assertions made in
Theorem~\ref{th-I3} that concern
$\rho_\alpha^{\log}$ and $\rho_{c,\alpha}^{\exp}$.
\begin{exa} We can apply Theorem~\ref{th-superpol-low1} when
\[
\rho(t)=\rho_{\alpha}^{\log}(t)=[\log(e+ t)]^\alpha,\qquad
\alpha>1.
\]
In this case we can take $\ell\simeq\rho_{\alpha-1}^{\log}$
and $\rho_1\simeq\rho_{\beta+1}^{\log}$ with $\beta>\alpha$. Then
$\widehat{\rho}_1\simeq\rho_\beta^{\log}$ and the conclusion is
that for any $\beta>\alpha$,
there are constants $C_2,c_\beta\in(0,\infty)$ such that,
for all $n$ large enough,
\[
- C_2n/ [\log n]^{\alpha-1} \le\log\Phi_{G,\rho_{\alpha}^{\log}}(n)
\le-c_\beta n/ [\log n]^{\beta}.
\]
\end{exa}
\begin{exa} Theorem~\ref{th-superpol-low1} gives a good result when
\[
\rho(t)=\rho_{c,\alpha}^{\exp}(t)=\exp\bigl(c[\log(1+
t)]^\alpha\bigr),\qquad
\alpha\in(0,1), c>0.
\]
Indeed, in this case we can obviously write
$\rho(t) = \log(e+t)\ell(t)$ with
$\ell\le\rho_{c,\alpha}^{\exp}$, and we can take
$\rho_1=\rho_{c_2,\alpha}^{\exp}$ for any fixed constant $c_2>c$.
The conclusion is
that there are constants $c_3,C_3$ such that,
for all $n$ large enough,
\[
- C_3n\exp(-c_3 [\log n]^{\alpha}) \le\log\Phi_{G,\rho_{c,\alpha
}^{\exp}}(n) \le-c_3n\exp(-C_3 [\log n]^{\alpha}).
\]
\end{exa}
\begin{theorem}
Assume that $G$ is a compactly generated locally compact
group and that there exist two continuous increasing functions
$\pi, \bar{\pi} $ such that, for all $n$ large enough,
\[
- n/\pi(n) \le \log\Phi_G(n)\le
-c n/\bar{\pi}(n).
\]
Assume that $\pi(t)\le t^{1-\varepsilon}$ for some $\varepsilon\in(0,1)$
Let $\rho(t)\simeq \log(e+t)\ell(t)$ where $\ell$ is a continuous increasing
slowly varying function at infinity.
Let $\rho_1$ be a slowly varying function such that
$\sum_1^\infty\frac{\rho(n)}{n\rho_1(n)}<\infty$ and set
$\widehat{\rho}
_1(t)=1/\int_t^\infty\frac{ds}{s\rho_1(s)}$.
Then there are $C_1,c_1(\rho_1)\in(0,\infty)$ such that,
for all $n$ large enough,
\[
-C_1 n/\ell(\pi(n^{\varepsilon/2}))
\le \log\Phi_{G,\rho}(n) \le- c(\rho_1) n/\widehat{\rho}_1(\bar{\pi}(n)).
\]
\end{theorem}
\begin{exa} Assume that $G= F\wr H$ where $F$ is a nontrivial
finite group and $H$ is polycylic with exponential volume growth. Then
$\Phi_{G}(n)\simeq\exp(-n/(\log n)^2)$. Hence, for any $\alpha>1$ and
$\beta>\alpha$, we obtain
\[
- C_1 n/ [\log(\log n)]^{\alpha-1}\le \Phi_{G,\rho^{\log}_\alpha}(n)
\le-c_\beta n/[\log(\log n)]^{\beta}
\]
for all $n$ large enough.
\end{exa}

\subsection{Volume-based upper bounds}
\label{sec-UpVol}

Let $G$ be a locally compact unimodular group equipped with its Haar
measure $\lambda$ (this group may well not be compactly generated).
Consider the problem of studying the decay of convolution powers of
probability measures of the form
%
\begin{equation}
\mu= \sum_1^\infty p_i \mu_i,\vadjust{\goodbreak}
\end{equation}
where $p_i\ge0$, $\sum_1^\infty p_i=1$ and
\[
\mu_i=\phi_i \,d\lambda,\qquad
\|\phi_i\|_\infty=\beta_i,\qquad \phi_i\ge0,\qquad \mu_i(G)=1.
\]
In words, $\mu$ is a convex linear
combination of the probability measures $\mu_i$, $i=1,2,\ldots,$ and
these measures are assumed to have bounded densities.
It was observed in~\cite{SCconv,Varconv} that interesting
upper bounds for convolution powers of such measures can sometimes
be obtained by elementary means. This is developed further below.

Set
\[
\sigma_k= \sum_{i> k}p_i,\qquad k=0,1,\ldots, \sigma_{-1}=+\infty,
\]
and
\[
b_k=\min_{i\le k}\{\beta_i\},\qquad k=1,2,\ldots, b_0=b_1,
\]
and consider the function $F$ on $(0,\infty)$ [this function depends
only on
$(\sigma_i)_0^\infty$ and $(b_i)_0^\infty$]
defined by
\[
F(s) = b_k \qquad\mbox{if } \sigma_{k}<s\le\sigma_{k-1}.
\]
The following result is quite versatile and surprisingly sharp when
applied to low moment measures.
\begin{pro}\label{prop4.8}
Referring to the notation introduced above and assuming that
$b_i\ra0$,
the density $\phi^{(n)}=d\mu^{(n)}/d\lambda$ of the $n$th convolution
power $\mu^{(n)}$
of $\mu$ satisfies
\[
\bigl\|\phi^{(n)}\bigr\|_\infty\le\int_0^\infty e^{-ns} \,dF(s)
=\sum_{i=1}^{\infty} e^{-n\sigma_i}
(b_i - b_{i+1}).
\]
\end{pro}
\begin{rem} One important class of examples is obtained by considering
a given increasing sequence of compact sets $B_i$ with $\bigcup_1^\infty
B_i=G$ and setting
\[
d\mu_i= d\lambda_{B_i}=\frac{1}{\lambda(B_i)}\mathbf1_{B_i}\,d\lambda.
\]
In this case, $b_i=\beta_i=1/\lambda(B_i)$.
\end{rem}
\begin{pf*}{Proof of Proposition~\ref{prop4.8}}
Write
\begin{eqnarray*}
\phi^{(n)}&=& \Biggl(\sum_{i=1}^\infty p_i\phi_i\Biggr)^{(n)}
= \sum_{k=1}^\infty\biggl(\biggl(\sum_{i\le k} p_i\phi_i\biggr)^{(n)}-\biggl(\sum_{i\le k-1}
p_i\phi_i\biggr)^{(n)}\biggr)\\
&=& \sum_{k\ge1} \biggl(\biggl(\sum_{i_j\le k} p_{i_1}\phi_{i_1}*\cdots*
p_{i_n}\phi_{i_n}\biggr)-
\biggl(\sum_{i_j\le k-1} p_{i_1}\phi_{i_1}*\cdots*p_{i_n}\phi_{i_n}\biggr)\biggr)\\
&=& \sum_{k\ge1} \biggl(\sum_{\max\{i_1,\ldots,i_n\}= k} p_{i_1}\phi
_{i_1}*\cdots*p_{i_n}\phi_{i_n}\biggr).
\end{eqnarray*}
Next we use Minkowski inequality and the estimate
\[
\| f_1*\cdots*f_n\|
_\infty\le\min\{\|f_i\|_\infty\}
\]
for functions $f_i$ with $L^1$-norm at most $1$. This estimate holds on
$G$ because we assume
unimodularity of $G$.
It yields
\begin{eqnarray*}
\bigl\|\phi^{(n)}\bigr\|_\infty&\le& \sum_{k\ge1} b_k\sum_{\max\{i_1,\ldots,i_n\}
= k} p_{i_1}\cdots p_{i_n}\\
&=& \sum_{k\ge1} b_k[ (1-\sigma_{k})^n- (1-\sigma_{k-1})^n]\\
&= & \sum_{k\ge1} (1-\sigma_k)^n [b_k-b_{k+1}]\\
&\le& \sum_{k\ge1} e^{-n\sigma_k} [b_k-b_{k+1}].
\end{eqnarray*}
\upqed\end{pf*}

We now give some simple applications when $G$ is locally compact,
compactly generated and unimodular (we assume that $G$ is noncompact).
Fix a symmetric open set $U$
that contains a generating compact neighborhood of the identity
element, and set
$|x|=\inf\{n\dvtx x\in U^n\}$, with $|e|=0$. Thus \mbox{$|\cdot| $} induces a familiar
word distance on $G$ when $G$ is finitely generated.
Observe that we have $\lambda(U^{4n})\ge2\lambda(U^n)$. Indeed,
if $|z|=3n$ (such a $z$ does indeed exist!), then the sets $U^n$ and $z U^n$
are disjoint and contained in $U^{4n}$. We consider the probability densities
$\phi_i=\lambda(B_i)^{-1}\mathbf1_{B_i}$ with $B_i=U^{4^i}$ and
set $b_i=\lambda(B_i)^{-1}$. Since $\lambda(B_i)\ge2\lambda(B_{i-1})$,
we have $b_i\ge b_i-b_{i+1}\ge b_i/2$. Set $\phi=\sum_1^\infty p_i\phi_i$
with $\sum_1^\infty p_i=1$ and $\sigma_k=\sum_{i>k}p_i$.
Fix a nondecreasing function $\rho\dvtx(0,\infty)\ra(0,\infty)$, and set
$\rho_G=\rho(|\cdot|)$. With this notation, we have
%
\begin{equation}\label{cvol}
\forall n\qquad \phi^{(2n)}(e)\le\sum_{k\ge1} e^{-2n\sigma_k} b_k
\end{equation}
and
%
\begin{equation}\label{cmom}
\int_G \rho_G \phi \,d\lambda= \int_G \sum_{k\ge1} p_k \rho_G \phi_k
\,d\lambda
\le \sum_{k\ge1}\rho(4^k) p_{k}.
\end{equation}
Further, we also have
\[
s\int_{\{\rho_G\ge s\}} \phi \,d\lambda\le s\sum_{\rho(4^{k-1})\ge s }
p_k .
\]
Hence, assuming that $\rho$ is a doubling function,
we have
%
\begin{equation}\label{wmom}
W(\rho,\phi \,d\lambda)=\sup_{s>0}\biggl\{
s\int_{\{\rho_G\ge s\}} \phi \,d\lambda\biggr\}
\le C(\rho)
\sup_{k}\{ \rho(4^k)\sigma_k\}.
\end{equation}
These two estimates allow us to derive upper bounds on $\Phi_{G,\rho_G}$
in terms of the volume growth of the group $G$. Indeed, for a given
$\rho$,
(\ref{cmom}) tells us how to pick $(p_i)_1^\infty$
so that $d\mu=\phi \,d\lambda$ satisfies $\mu(\rho_G)<\infty$.
For this choice of $(p_i)_1^\infty$,~(\ref{cvol}) yields an upper bound
on $\phi^{(2n)}(e)$ [hence on $\Phi_{G,\rho_G}(n)$] in terms of
(a lower bound on) the volume
growth which determines the sequence $(b_i)_1^\infty$.
This approach yields an alternative proof of Theorem~\ref{th-Phi-Pol1}
(polynomial volume growth) as well as new results in
the super-polynomial volume case.
\begin{pf*}{Alternative proof of Theorem~\ref{th-Phi-Pol1}}
We give the details only for $\widetilde{\Phi}_{G,\rho}$. The proofs
concerning $\Phi_{G,\rho}$ are similar.
Recall that Theorem~\ref{th-Phi-Pol1} deals with the case when
$V(n)\simeq n^D$,
and $\rho$ is comparable to either (a) a
regularly varying function with positive index $\alpha\in(0,2)$
or (b) a slowly varying of the form $\rho(t)\simeq
1/\int_t^\infty\frac{ds}{s\ell(s)}$
with $\ell$ positive and slowly varying. Further, in case~(b), we
assume that
$\log\rho^{-1} (t)\simeq t^\gamma\omega(t)^{1+\gamma}$ for some
$\gamma\in[0,\infty)$ and positive slowly varying function $\omega$.

In case (a), set $p_i= c \rho(4^i)^{-1}$. In case (b) set
$p_i=c\ell(4^i)^{-1}$. Then it is easy to check that
\[
\sigma_i = \sum_{k> i}p_k \simeq
\rho(4^i)^{-1}.
\]
Using~(\ref{wmom}), this implies that $\phi=\sum_1^\infty p_i\phi_i$ satisfies
the moment condition
\[
W(\rho,\phi \,d\lambda)\le\sup_i\{\rho(4^i)\sigma_i\} <\infty.
\]
Further
%
\begin{equation}\label{In0}
\phi^{(2n)}(e) \le
C_1\sum_i e^{- c_1 n /\rho(4^{i})} 4^{-i D} .
\end{equation}
In case (a) where $\rho(t)\simeq(1+t)^{2\alpha}\ell(s^2)^\alpha$ with
$\alpha\in(0,1)$ and $\ell$ positive and slowly varying, observe that
$\rho^{-1}(1/u)\simeq u^{-1/2\alpha}\ell^{\#}(1/u^{1/\alpha})^{1/2}$
for small
$u$ and write
\begin{eqnarray*}
\phi^{(2n)}(e) &\le&
C_1\sum_i e^{- c_1 n /\rho(4^{i})} 4^{-i D} \leq
C_2\int_1^\infty e^{-c_1 n/\rho(s)}\,\frac{ds}{s^{1+D}}\\[-1pt]
&\le& C_3\int_0^1 e^{-c_1 n
u}\biggl(\frac{1}{\rho^{-1}(1/u)}\biggr)^{D} \,\frac
{du}{u}\\[-1pt]
&\le& C_4 \int_0^\infty e^{-c_1 n u} \biggl(\frac{u^{1/\alpha}}{\ell^{\#
}(1/u^{1/\alpha})}\biggr)^{D/2}\,\frac{du}{u}
\simeq C_5 [n^{1/\alpha}\ell^{\#}(n^{1/\alpha})]^{-D/2}.
\end{eqnarray*}
This yields the desired result, namely,
\[
\widetilde{\Phi}_{G,\rho}(n)\le C[n^{1/\alpha}\ell^{\#}(n^{1/\alpha})]^{-D/2}
\]
for case (a).

In case (b), write
\begin{eqnarray*}
\phi^{(2n)}(e) &\le&
C_1\sum_i e^{- c_1 n /\rho(4^{i})-(D/2)\log 4^{-i }} 2^{-Di}\\
&\le& C_2\exp\Bigl( -c_2\inf_{s>0} \{n/\rho(s) +\log(e+ s)\}\Bigr).
\end{eqnarray*}
Using the hypothesis concerning $\rho^{-1}$, observe that
\begin{eqnarray*}\inf_{s>0} \{n/\rho(s) +\log(e+ s)\}&=&
\inf_{s>0} \bigl\{ns +\log\bigl(e+ \rho^{-1}(1/s)\bigr)\bigr\}\\
&\simeq&
\inf_{s>0} \{ns + s^{-\gamma} \omega(1/s)^{1+\gamma}\}\\
&\simeq& n^{\gamma/(1+\gamma)}/\omega^{\#}\bigl(n^{1/(1+\gamma)}\bigr).
\end{eqnarray*}
As stated in Theorem~\ref{th-Phi-Pol1}(2) and under the hypotheses of case (b),
this yields
\[
\widetilde{\Phi}_{G,\rho}(n) \le
C\exp\bigl(- c n^{\gamma/(1+\gamma)}/\omega^{\#}\bigl(n^{1/(1+\gamma)}\bigr)\bigr)
\]
as desired.
\end{pf*}
\begin{theorem}[(The super polynomial case)] \label{th-Vol2}
Assume that $\lambda(U^n)\ge\exp(cn^\theta)$
for some $c,\theta>0$. Then:
\begin{longlist}[(4)]
\item[(1)] Fix $\alpha\in(0,1)$, a positive slowly
varying function $\ell$ at infinity, and set
$\rho(s)=[(1+s^2) \ell(1+s^2)]^\alpha$. Then there exists
$C\in(0,\infty)$ such that, for all $n$ large enough,
\[
\widetilde{\Phi}_{G,\rho_\alpha}(n) \le C \exp\bigl(-c
n^{\theta/(\theta+2\alpha)}/
\ell_\bullet^{\#}\bigl(n^{2/(\theta+2\alpha)}\bigr)^{\alpha}\bigr),
\]
where $\ell_\bullet=[\ell^{\#}]^{\theta/(\theta+2\alpha)}$,
and $\ell^{\#}_\bullet$ is its de Bruijn conjugate.
\item[(2)] Fix $\alpha\in(0,2)$. For all $\beta>\alpha$,
there are constants $C_\beta,c_\beta>0$
such that, for all $n$ large enough, $\Phi_{G,\rho_\alpha}(n) \le
C_\beta\exp(-c_\beta
n^{\theta/(\theta+\beta)})$.\vspace*{1pt}
\item[(3)] For any fixed $\alpha>0$ we have
$\widetilde{\Phi}_{G,\rho^{\log}_\alpha}(n) \le C
\exp(-c n/[\log n]^\alpha)$.
Further, for all $\beta>\alpha$,
there are constants $C_\beta,c_\beta>0$
such that, for all $n$ large enough,
$\Phi_{G,\rho^{\log}_\alpha}(n) \le C_\beta
\exp(-c_\beta n/[\log n]^\beta)$.\vspace*{1pt}
\item[(4)] For any fixed $\alpha\in(0,1)$ and $c>0$,
there is a constant $C_1>0$
such that, for all $n$ large enough, $\Phi_{G,\rho^{\exp}_{c,\alpha
}}(n) \le C_1
\exp(- n/ \exp(C_1[\log n]^\alpha)))$.
\end{longlist}
\end{theorem}
\begin{pf} We prove statement (1).
The variations needed for the other statements are straightforward.
We have
$b_i \le\exp(-c 4^{i\theta})$.
Fulfilling the desired moment conditions forces the
choice\vadjust{\goodbreak} of the sequence $(p_i)_1^\infty$.
For instance, in the first case, we take $p_i =c \rho( 4^i)^{-1}$
so that $\sigma_i\simeq\rho(4^i)^{-1}$. Hence
\begin{eqnarray*}
\phi^{(2n)}(e)
&\le& C_1\sum_i e^{- c_1( n /\rho(4^{i}) + 4^{i\theta})} \\
&\le&C_2 \exp\Bigl(- c_2 \inf_{s>0}\{ n /\rho(s) + s^\theta
\}\Bigr).
\end{eqnarray*}
Write
\[
\inf_{s>0}\{ n /\rho(s) + s^\theta
\}=\inf_{s>0}\{ n s + \rho^{-1}(1/s)^\theta\}.
\]
A good\vspace*{1pt} approximation of the infimum is obtained by picking
$s=s_n$ such that $n= (1/s_n)\rho^{-1}(1/s_n)^\theta$.
At infinity, $\rho^{-1}(t) = t^{1/2\alpha}\ell^{\#}(t^{1/\alpha})^{1/2}$
and thus, at $0$,
\[
(1/t)[\rho^{-1}(1/t)]^\theta=
t^{-(2\alpha+\theta)/2\alpha}\ell^{\#}(1/t^{1/\alpha})^{\theta/2}.
\]
Setting
$\ell_\bullet=[\ell^{\#}]^{\theta/(\theta+2\alpha)}$, we have
$s_n\simeq n^{-2\alpha/(2\alpha+\theta)}[\ell_\bullet^{\#
}(n^{2/(2\alpha+\theta)})]^{-\alpha}$.
Finally,
\[
\phi^{(2n)}(e)\le C \exp\bigl(- c_3 n^{\theta/(\theta+2\alpha)}/
\ell_\bullet^{\#}\bigl(n^{2/(\theta+2\alpha)}\bigr)^{\alpha}\bigr).
\]
\upqed\end{pf}
\begin{rem}
Note that
the hypotheses in Theorem~\ref{th-Vol2} and in
Theorem~\ref{th-superpol} are notably different.
Theorem~\ref{th-superpol}
is based on hypotheses regarding the behavior of $\phi_G$ whereas
Theorem~\ref{th-Vol2} assumes $V(n) \ge\exp(c n ^\theta)$.
We note that the hypothesis
$V(n) \ge\exp(c n ^\theta)$ implies
$\Phi_G(n)\le\exp(-c n^{\theta/(2+\theta)})$~\cite{VSCC}.
If $V(n)\ge\exp(c n^\theta)$ and $\Phi_G(n)\ge\exp(-Cn^{\theta
/(2+\theta)})$,
then the upper bound of Theorem~\ref{th-Vol2}(1) matches
precisely the lower bound of Theorem~\ref{th-Phi>2}.
\end{rem}

The next theorem treats the case of groups
that have exponential volume growth (i.e., $\theta=1$) and such that
$\Phi_G(n)\simeq\exp(- n^{1/3})$. (This is the case, e.g., if
$G$ is polycyclic with exponential volume growth.)
This result contains the part of Theorem~\ref{th-I3} dealing
with $\rho_\alpha$, $\alpha\in(0,2)$.
\begin{theorem} Assume that $G$ has exponential volume growth and satisfies
$\phi_G(n)\simeq\exp(-n^{1/3})$.
\begin{longlist}[(2)]
\item[(1)] Fix $\alpha\in(0,1)$, a positive slowly
varying function $\ell$ at infinity, and set
$\rho(s)=[(1+s^2) \ell(1+s^2)]^\alpha$. Then we have
\[
\widetilde{\Phi}_{G,\rho}(n)\simeq \exp\bigl(-
n^{1/(1+2\alpha)}/\ell_\bullet^{\#}
\bigl(n^{2/(1+2\alpha)}\bigr)^{\alpha}\bigr),
\]
where $\ell_\bullet=[\ell^{\#}]^{1/(1+2\alpha)}$, and
$\ell^{\#}_\bullet$ is its de Bruijn conjugate.
\item[(2)] Fix $\alpha\in(0,2)$.
For all $\beta>\alpha$,
there are constants $C_\beta,c_\beta>0$
such that, for all $n$ large enough,
\[
C_\beta\exp\bigl(-c
n^{1/(1+\alpha)}\bigr)
\le \Phi_{G,\rho_\alpha}(n) \le C_\beta\exp\bigl(-c_\beta
n^{1/(1+\beta)}\bigr).
\]
\end{longlist}
\end{theorem}

\section{The case of the wreath product $(\mathbb Z/2\mathbb Z)\wr
\mathbb Z^d$}\label{sec-Wreath}

The wreath product construction provides important examples of groups
whose behavior differs from linear groups. The simplest family
of wreath products is\vadjust{\goodbreak} $(\mathbb Z/2\mathbb Z)\wr\mathbb Z^d$.
An element of this group is a pair $(\eta,k)$ with
$\eta\in\bigoplus_{i\in\mathbb Z^d} (\mathbb Z/2\mathbb Z)_i$ (algebraic sum)
and $k\in\mathbb Z^d$. In the popular lamplighter interpretation,
$k$ is the position of the lamplighter, and $\eta=(\eta_i)_{i\in
\mathbb Z^d}$
is a configuration of lamps that can be on ($\eta_i=1$) or off
($\eta_i=0$). Only finitely many lamps can be on. The product is given by
$(\eta,k)(\eta',k')= (\eta'',k'')$ where
$k''=k+k'$ (addition in $\mathbb Z^d$) and $\eta''_i= \eta_i+\eta'_{i-k}$
(addition in
$\mathbb Z/2\mathbb Z$). In other words, $(\mathbb Z/2\mathbb Z)\wr
\mathbb Z^d$
is the semidirect product of $\bigoplus_{i\in\mathbb Z^d} (\mathbb
Z/2\mathbb Z)_i$ by $\mathbb Z^d$ where the action of $\mathbb Z^d$ on
$\bigoplus_{i\in\mathbb Z^d} (\mathbb Z/2\mathbb Z)_i$ is by translation
of the indices. These groups have exponential volume growth.

The aim of this section is to prove the following theorem.
\begin{theorem} \label{th-wr}
For any integer $d\ge1$ and $\alpha\in(0,2)$, we have
\[
\widetilde{\Phi}_{(\mathbb Z/2\mathbb Z)\wr\mathbb Z^d,\rho_\alpha
}(n)\simeq
\exp\bigl(- n^{d/(d+\alpha)}\bigr).
\]
Further, for any $\beta>\alpha$, there are constants
$c,C,c_\beta,C_\beta\in(0,\infty)$ such that, for all $n$ large enough,
\[
c\exp\bigl(-C n^{d/(d+\alpha)}\bigr) \le
\Phi_{(\mathbb Z/2\mathbb Z)\wr\mathbb Z^d,\rho_\alpha}(n)\le
C_\beta\exp\bigl(-c_\beta n^{d/(d+\beta)}\bigr).
\]
\end{theorem}

We shall see in the proof given below that
the lower bounds stated in this theorem follow from Theorem~\ref{th-Phi>2}.
The interesting part are the upper bounds. These upper bounds are interesting
because they do not follow from the results in Sections~\ref{sec-Sub} and
\ref{sec-UpVol}.
\begin{pf*}{Proof of Theorem~\ref{th-wr}}
We can identify $\mathbb Z^d$ as a subgroup of $(\mathbb Z/2\mathbb
Z)\wr\mathbb Z^d$ in an obvious way, and we can also identify $\mathbb
Z/2\mathbb Z$ with $(\mathbb Z/2\mathbb Z)_0$ in
$\bigoplus_{i\in\mathbb Z^d}(\mathbb Z/2\mathbb Z)_i \subset(\mathbb
Z/2\mathbb Z)\wr\mathbb Z^d$. Hence, any probability measure on
$\mathbb Z/2\mathbb Z$ or on $\mathbb Z^d$ can be interpreted as a
measure on $(\mathbb Z/2\mathbb Z)\wr\mathbb Z^d$. Following the
notation used in~\cite{PSCwp}, if $\nu$ is a measure supported on
$(\mathbb Z/2\mathbb Z)_0$, and $\mu$ a measure supported on $\mathbb
Z^d$, we set $q=\nu*\mu*\nu$ in $(\mathbb Z/2\mathbb Z)\wr\mathbb Z^d$.
In~\cite{PSCwp}, it is observed that a famous large deviation theorem,
due to Donsker and Varadhan~\cite{DV} and concerning the range of
certain random walks on~$\mathbb Z^d$, implies that
\[
q^{(2n)}(e)\simeq\exp\bigl(-n^{d/(d+2)}\bigr),
\]
when $\nu$ is the uniform measure on $\mathbb Z/2\mathbb Z$, and $\mu$
is any
symmetric measure on $\mathbb Z^d$ with finite generating support.
By~\cite{PSCstab}, this implies that
%
\begin{equation}\label{Phi-wr}
\Phi_{(\mathbb Z/2\mathbb Z)\wr\mathbb
Z^d}(n)\simeq\exp\bigl(-n^{d/(d+2)}\bigr).
\end{equation}

Here, we are interested in determining the behavior of
\[
\widetilde{\Phi}_{(\mathbb Z/2\mathbb Z)\wr\mathbb Z^d,\rho_\alpha}
\quad\mbox{and}\quad
\Phi_{(\mathbb Z/2\mathbb Z)\wr\mathbb Z^d,\rho_\alpha},\qquad
\alpha\in(0,2).
\]
First, consider how the results obtained so far in this paper apply in
this case.
Theorem~\ref{th-Phi>2} readily gives the lower bound
%
\begin{equation}\label{wr-low}
\Phi_{(\mathbb Z/2\mathbb Z)\wr
\mathbb Z^d,\rho_\alpha}(n)\ge
\widetilde{\Phi}_{(\mathbb Z/2\mathbb Z)\wr
\mathbb Z^d,\rho_\alpha}(n)\ge
\exp\bigl(-C(d,\alpha)n^{d/(d+\alpha)}\bigr),\vadjust{\goodbreak}
\end{equation}
because if $\gamma=d/(d+2)$, then
$\gamma_{\alpha/2} :=\gamma/(\gamma+(\alpha/2)(1-\gamma))=
d/(d+\alpha)$.
We are faced with the problem of deciding whether or not this is sharp.
Can we find measures with finite $\rho_\alpha$-moment and whose convolution
powers decay as rapidly as permitted by this lower bound?

For this purpose, we have so far discussed two methods: (a) the use of
subordination as developed in~\cite{BSCsubord} and (b) direct computation
based on volume estimates (see Theorem~\ref{th-Vol2}).

The direct computation of Theorem~\ref{th-Vol2} provides the upper bounds
%
\begin{equation} \label{wr-upV1}
\widetilde{\Phi}_{(\mathbb Z/2\mathbb Z)\wr\mathbb Z^d,\rho_\alpha
}(n)\le
\exp\bigl(-C
n^{1/(1+\alpha)}\bigr)
\end{equation}
and
%
\begin{equation} \label{wr-upV}
\Phi_{(\mathbb Z/2\mathbb Z)\wr\mathbb Z^d,\rho_\alpha}(n)\le\exp
\bigl(-C(\beta)
n^{1/(1+\beta)}\bigr),\qquad \beta>\alpha.
\end{equation}
When $d=1$ (and only in this case), these upper bounds show that
the lower bounds stated in~(\ref{wr-low}) are essentially sharp.
In particular, we get
\[
\mathrm{exp}\mbox{-}\mathrm{pow}\bigl((\mathbb Z/2\mathbb Z)\wr\mathbb Z,\rho_\alpha\bigr)=
1/(1+\alpha),\qquad \alpha\in(0,2).
\]
For $d\ge2$,~(\ref{wr-upV1}) and~(\ref{wr-upV}) fail to match
(\ref{wr-low}) for a good reason: Theorem~\ref{th-Vol2} is based solely
on a volume hypothesis and thus cannot provide more subtle information
that is based on the particular structure of these wreath products.

The subordination technique of~\cite{BSCsubord} fails to give good
upper bounds for a different reason related to the fact that, for
simple random walks on wreath products such as $(\mathbb Z/2\mathbb
Z)\wr\mathbb Z^d$ with $d>1$, the rate of escape to infinity is much
faster than~$\sqrt{n}$. See the discussion in~\cite{BSCsubord}.

Thus, the two techniques used earlier in this paper
to provide upper bounds on $\Phi_{G,\rho_\alpha}$ and
$\widetilde{\Phi}_{G,\rho_\alpha}$
both fail to match the lower bound~(\ref{wr-low}) when $d\ge2$.
The following argument shows that~(\ref{wr-low}) is
sharp nonetheless.
For each $\alpha\in(0,2)$ let $\mu_{\alpha}$ be the probability
measure on $\mathbb Z^d$ given by
\[
\mu_\alpha(k)= \frac{c(d,\alpha)}{(1+\|k\|^2)^{(d+\alpha)/2}},\qquad
k\in\mathbb Z^d, \|k\|^2=\sum_1^d k_i^2.
\]
The theorem of Donsker and Varadhan (\cite{DV}, Theorem 1) implies
that, for any fixed $s$ and $n$ large enough,
\[
E(e^{ -s D^{\#}_n}) \simeq\exp\bigl(- n^{d/(d+\alpha)}\bigr).
\]
Here $D^{\#}_n$ is the number of visited sites up to time $n$ for the
random walk on $\mathbb Z^d$ driven by $\mu_\alpha$.

For any fixed $\beta\in(0,2)$, this, together with
\cite{PSCwp}, Theorem 3.1,
implies that the
measure $q_\beta=\nu*\mu_\beta*\nu$ on $(\mathbb Z/2\mathbb Z)\wr
\mathbb Z^d$
satisfies
\[
q^{(2n)}(e)\simeq\exp\bigl(- n^{d/(d+\beta)}\bigr).
\]
It is plain that the measure $q_\alpha$ has finite
weak-$\rho_\alpha$-moment $W(\rho_\alpha,q_\alpha)<\infty$,
and that $q_\beta$ has finite\vadjust{\goodbreak} $\rho_\alpha$-moment
if and only if $\beta>\alpha$. To check this, notice that $q_\alpha$ is
almost entirely concentrated on $\mathbb Z^d$ inside
$(\mathbb Z/2\mathbb Z)\wr\mathbb Z^d$.
Thus these measures provide witnesses
to the fact
that
\[
\widetilde{\Phi}_{(\mathbb Z/2\mathbb Z)\wr\mathbb Z^d,\rho_\alpha}
(n)\le
C_1\exp\bigl(-c_1 n^{d/(d+\alpha)}\bigr)
\]
and that,
for each $\beta>\alpha$, $\alpha\in(0,2)$,
\[
\Phi_{(\mathbb Z/2\mathbb Z)\wr\mathbb Z^d,\rho_\alpha} (n)\le
C_\beta\exp\bigl(-c_\beta n^{d/(d+\beta)}\bigr).
\]
These are the desired upper bounds.
\end{pf*}

In particular it follows that
\[
\mathrm{exp}\mbox{-}\mathrm{pow}\bigl((\mathbb Z/2\mathbb Z)\wr\mathbb
Z^d,\rho_\alpha\bigr)= \frac{d}{d+\alpha},\qquad d=1,2,\ldots,
\alpha\in(0,2).
\]
It is interesting to note that the optimal measure $q_\alpha$
that we have exhibited above is spread out only in a very small part of
the group, that is, in the directions of the lamplighter moves $\mathbb Z^d$.

\begin{appendix}\label{app}

\section*{Appendix: Ultracontractivity, functional calculus
and von Neumann trace}

\subsection{Spectral theory} Let $T$ be a self-adjoint operator
acting on a Hilbert space $H$. We denote by $E^T_{\mathcal I}$
its spectral projector
associated with the open set $\mathcal I\subset\mathbb R$,
and by $E^T_s=E^T_{(-\infty,s)}$
the associated (left-continuous) spectral resolution of $T$ so that
\[
T=\int_{-\infty}^{+\infty} s\,dE^T_s.
\]

In the cases\vspace*{1pt} of interest to us, $T$ is actually a bounded operator so
that $E^T_{(a,b)}=0$ if $\max\{a,-b\}$ is larger than $\|T\|$.
For any continuous function $m\dvtx\mathbb R\ra\mathbb C$, the operator
$m(T)$ with domain
$D_f=\{u\in H\dvtx\int|m(s)|^2\,d\langle E^T_su, u\rangle<\infty\}$
is defined by
\[
m(T)= \int_{-\infty}^{+\infty}m(s) \,dE^T_s,
\]
where this integral is obtained as the strong limit of finite Riemann sums.
Further, note that if $m$ is real valued, then $m(T)$ is self-adjoint and
\[
E^{m(T)}_{(a,b)}= E^T_{m^{-1} (a,b)}.
\]

\subsection{The von Neumann algebra $V(G)$}
We will make fundamental use of the notion of von Neumann trace for
certain operators in the von Neumann algebra $V(G)$
generated by the right translations $r_g\dvtx f\mapsto f(\cdot g)$
acting on $L^2(G)$. By construction, $V(G)$ is equipped with a faithful
semifinite normal trace $\tau$ defined as follows.
Let $S$ be a nonnegative Hermitian
element in $V(G)$ [i.e., a self-adjoint\vadjust{\goodbreak} element satisfying
$\langle Su,u\rangle\ge0$ for every $u\in L^2(G)$]. If $S^{1/2}=R_a$
for some $a\in L^2(G)$, set $\tau(S)=\|a\|_2^2$. Otherwise,
set $\tau(S)=+\infty$. See~\cite{Dix2}, page 97. Since $S^{1/2}=R_a$,
$R_a$ is self-adjoint.
This is equivalent to say that the function $a\in L^2(G)$
satisfies $a=\check{a}$ [where $\check{a}(x)=\bar{a}(x^{-1})$, $x\in G$].
Hence
\[
\tau(S)= \int_G|a|^2\,d\lambda= a*a(e).
\]
Note that, as the convolution of two functions in $L^2(G)$,
the function $a*a$ is bounded and continuous
(i.e., admits a continuous representative) and that $S$ acts on
$\phi\in\mathcal C_c(G)$ by $S\phi= \phi*[a*a]$.

Let $ S,T$ be two Hermitian nonnegative elements in $V(G)$
such that $S\le T$. Then $\tau(S)\le\tau(T)$. In particular,
if $T$ has finite trace and spectral decomposition
\[
T=\int_0^\infty s\,dE_s^T,
\]
then $E^T_{(s,+\infty)}$ is in $V(G)$ and has finite trace for all
$s>0$ since
$sE^T_{(s,+\infty)}\le T$. Note that, in general (i.e., when $G$
is not countable), $E^T_\infty=I$ does not have finite trace.

If $T$ is Hermitian of the form
$T=R_{a*\check{a}}$,
then
\[
\tau(T)=a*\check{a}(e)= \int_{0}^{+\infty} s\,
d\bigl[-\tau\bigl(E^T_{(s,\infty)}\bigr)\bigr]=
\int_{0}^{+\infty} \tau\bigl(E^T_{(s,\infty)}\bigr)\,ds.
\]
This follows from the well-known properties of spectral resolutions
and the fact that $\tau$ is a normal trace [this means that $\tau$ has the
property that, for any positive Hermitian $T$ and any increasing
filtering set $\mathcal F$
of positive Hermitian elements with supremum $T$,
$\sup_{\mathcal F}\tau(S)=\tau(T)$].

Strictly speaking, the trace $\tau$
is defined only on nonnegative Hermitian elements. However, the set of
Hermitian nonnegative elements with finite trace is the positive part of
a two-sided ideal $\mathfrak m$ of $V(G)$, and
there is a unique linear form defined on this two-sided ideal
which coincides with the trace on nonnegative Hermitian elements.
Abusing notation, we denote this extension by
$\tau\dvtx\mathfrak m\ra\mathbb R$. If
$a,b\in L^2(G)$ and $R_a,R_b\in V(G)$, then $R_aR_b \in\mathfrak m$ and
$\tau( R_a*R_b)= b*a(e)$. See~\cite{Dix2}, Theorem~1, page 97. In particular,
if $\phi=\check{\phi}\in L^1(G)\cap L^2(G)$ and
$T=R_\phi$, then, for any $n=2,3,\ldots,$
$T^n=R_{\phi^{(n)}}$ has finite trace and
%
\setcounter{equation}{0}
\begin{equation}\label{trace=}
\phi^{(n)}(e)=\tau(T^n).
\end{equation}

\subsection{Ultracontractivity}

Let $\phi=\check{\phi}\in L^1(G)\cap L^2(G)$ be a symmetric
probability density.
Let $ T=R_{\phi}\dvtx f\mapsto f*\phi$ be the operator of right convolution by
$\phi$ acting on $L^2(G)$. This is an
Hermitian element of $V(G)$ with norm at most $1$.
Its powers $T^n$, $n\ge2$, are of finite trace and the function
\[
n\mapsto\tau(T^n)=\phi^{(n)}(e)
\]
is of interest to us because it quantifies the ultracontractivity of the
operators $T^{2n}$, $n\ge1$. Indeed, we have
\[
\sup_{\|f\|_1\le1}\{\| T^{2n}f\|_\infty\}=\|T^{2n}\|_{1\ra\infty}=
\phi^{(2n)}(e).
\]
We assume throughout that $\phi^{(2n)}(e)\ra0$, which simply
means that $\phi$ is not supported on a compact subgroup of $G$.
As a consequence $\|T^nf\|_\infty\ra0$ for any $f\in L^2(G)$. In particular,
there are no nontrivial functions in $L^2(G)$ such that $Tf=\pm f$.

Let $E^T_s,E^{I-T}_s$, $s\in\mathbb R$,
be the left-continuous spectral resolutions
of $T$ and $I-T$ and note that
\[
T^n=\int_0^{2}(1-s)^n\,d E^{I-T}_s,\qquad E^T_{(1-b,1-a)}=E^{I-T}_{(a,b)},
\qquad 0\le a<b\le\infty,
\]
with $\lim_{s\searrow0}E^{I-T}_s=E^T_{[1,\infty)}=0$ because there are
no $L^2(G)$-solutions of \mbox{$Tu=u$}. Note also that the projection valued
measure $dE^{I-T}_s$ could have an atom at $s=1$ [corresponding to
$L^2(G)$-functions satisfying $Tu=0$] but that this atom is irrelevant
to the integral formula
$m(T)=\int_0^2 m(1-s)\,dE^{I-T}_s$ as long as $m$ is continuous and
satisfies $m(0)=0$.
Observe further that $(1-s)^2E^{I-T}_s\le T^2$ for $s\in[0,1]$
so that $E^{I-T}_s$
has finite trace for all $s\in[0,1)$. Similarly $E^{I-T}_{(s,2)}$ has
finite trace
for $s\in(1,2)$.

Using this fact we define
the nondecreasing, nonnegative functions $N_{\phi}\dvtx \break [0, 1)\ra
[0,+\infty)$
by
%
\begin{equation}\label{def-Nphi}
N_{\phi}(s)=\tau( E^{I-T}_s)=\tau\bigl(E^{T}_{(1-s,\infty)}\bigr),\qquad
s\in(0,1).
\end{equation}
The following lemma is proved in~\cite{BSCsubord}. It indicates that
the part of the spectrum of $T$ near $-1$ does not play a crucial role in
estimating $\phi^{(2n)}(e)$ [this uses the fact that $\phi^{(k)}(e)\ge0$].
\setcounter{theorem}{0}
\begin{lem}[(See, e.g.,~\cite{BSCsubord}, Proposition 3.1)] \label{lem-21}
Assume that $\phi$ is a symmetric
probability density in $L^2(G)$. Then
%
\begin{equation}
\int_0^1(1-s)^{2n}\,dN_\phi(s)\le\phi^{(2n)}(e)\le2\int_0^1(1-s)^{2(n-1)}\,
dN_\phi(s) .
\end{equation}
\end{lem}

Thanks to this Laplace transform type relation, the behavior of
$n\mapsto\phi^{(2n)}(e)$ as $n$ tends to infinity and the behavior of
$N_{\phi}(s)$ as $s$ tends to $0$ are related to each other.
The following statements are appropriate versions of classical results.
See~\cite{BPS,BSCsubord,BGT} for details.

For $\theta=0$ or $+\infty$, we let $\mathcal R_\alpha(\theta)$
be the set of regularly varying functions of index $\alpha$ at $\theta$.
If $\ell$ is a slowly varying function at infinity, we let $\ell^{\#}$ be
its de Bruijn conjugate. See~\cite{BGT}, Theorem 1.5.13. For simple
applications,
we observe that if $\ell(x)\sim\ell(x\ell(x))$ at infinity,
then $\ell^{\#}\sim1/\ell$. For instance, this applies to
$\ell(x)=(\log x)^\beta$, $\beta\in\mathbb R$. See~\cite{BGT},
Corollary 2.3.4.\vadjust{\goodbreak}
In the following result, $\varphi$ and $N$ are abstract functions but,
applications we have in mind, $\varphi(k)=\phi^{(2k)}(e)$ and $N=N_\phi$
as in Lemma~\ref{lem-21}.
\begin{pro} \label{pro-NPhi2}
Let the nondecreasing function $N\dvtx(0,1)\ra(0,+\infty)$
and nonincreasing function $\varphi\dvtx \{1,2,\ldots\}\ra(0,+\infty)$ be
related by
\[
\forall k>k_0\qquad c\int_0^1(1-s)^{k}\,dN(s)\le\varphi(k) \le C \int
_0^1(1-s)^{k-k_0}\,dN(s)
\]
for some $k_0,c,C\in(0,\infty)$.
\begin{longlist}[(3)]
\item[(1)] Fix $\alpha>0$, and let $\ell$ be a slowly varying function
at infinity.
There exists a $c_1\in(0,1)$ such that
$\varphi(k)k^{\alpha}\ell(k)\ge c_1 $
[resp., $\varphi(k)k^{\alpha}\ell(k)\le c_1$]
for all $k$ large enough
if and only if
there exists a constant $c_2\in(0,1)$ such that
$N(s)s^{-\alpha}\ell(1/s) \ge c_2$ [resp., $N(s)s^{-\alpha}\ell(1/s)
\le c_2$]
for all $s>0$ small enough.

\item[(2)] Fix $\alpha\in(0,1)$, and let $\ell$ be a slowly varying
function at
infinity.
There exists a constant $c_1\in(0,1)$ such that
\[
[-\log\varphi(k)] [k^\alpha/\ell(k^{1-\alpha})]^{-1}\ge c_1
\qquad(\mbox
{resp., $\le$}
c_1)\qquad
\mbox{ for large enough } k,
\]
if and only if there exists a constant $c_2\in(0,1)$ such that
\begin{eqnarray}
[-\log N(s)] [s^{\alpha} \ell^{\#}(1/s)]^{1/(1-\alpha)}\ge c_2
\qquad(\mbox{resp., $\le$} c_2)\nonumber\\
&&\eqntext{\mbox{for small enough } s>0.}
\end{eqnarray}
\item[(3)] Let $M$, $\pi$ and $t\mapsto t/\pi(t)$ be continuous
increasing functions
on $(0,\infty)$ which tend to infinity at infinity and such that
%
\begin{equation} \label{rel-Mpi*}
\pi^{-1}(t)\simeq t M(t) \qquad\mbox{at infinity}.
\end{equation}
The following two properties are equivalent:
\begin{enumerate}[(a)]
\item[(a)] there exists $c_1\in(0,\infty)$ such that
\[
-\log N(s)\ge c_1 M(c_1/s)\qquad
\mbox{[resp., $-\log N(s)\le c_1 M(c_1/s)$]}
\]
for
all $s$ small enough;
\item[(b)] there exists $c_2\in(0,\infty)$ such that
\[
-\log\varphi(n) \ge c_2 n/\pi(n/c_2)\qquad \mbox{[resp., $-\log\varphi(n) \le
c_2 n/\pi(n/c_2)$]}
\]
for all $n$ large enough.
\end{enumerate}
\end{longlist}
\end{pro}
\begin{exa} The reason behind considering these elaborate
statements is the nature of the known results concerning $\phi^{(2n)}(e)$
when $\phi$ is symmetric compactly supported. Here is a
small selection of specific examples of interest.
\begin{longlist}[(3)]
\item[(1)] The properties
\[
\phi^{(2n)}(e)\simeq n^{-D/2} \qquad\mbox{at infinity}
\]
and
\[
N_\phi(s) \simeq s^{D/2} \qquad\mbox{at zero}\vadjust{\goodbreak}
\]
are equivalent. These properties
hold when $\phi$ is compactly supported, and $G$ has polynomial volume
growth of degree $D$.
\item[(2)] The properties
\[
\phi^{(2n)}(e)\simeq \exp(-n^{1/3})\qquad
\mbox{at infinity}
\]
and
\[
N_\phi(s)
\simeq\exp(-1/s^{1/2}) \qquad\mbox{at zero}
\]
are equivalent.
These properties hold whenever $\phi$ is compactly supported with generating
support and
$G$ is virtually polycyclic with exponential volume growth.
\item[(3)] The properties
\[
\phi^{(2n)}(e)\simeq
\exp\bigl(-n ^{d/(d+2)}[\log n]^{2/(d+2)}\bigr) \qquad\mbox{at infinity}
\]
and
\[
N_\phi(s)\simeq\exp(- s^{-d/2}[\log1/s]) \qquad\mbox{at zero}
\]
are equivalent.
They hold, for instance, when $G=\mathbb Z \wr\mathbb Z^d $ (the
lamplighter group with
street map $\mathbb Z^d$ and lamps in $\mathbb Z$).
\end{longlist}
See~\cite{Alex,BPS,Erschleriso,ErschlerPTRF,PSCwp,VSCC}.
Remarkably enough, the first two types of behaviors
are the only possibilities for unimodular amenable Lie groups and for
finitely generated amenable discrete subgroups of Lie groups; see, for example,
\cite{SCsurveydiff} and the references therein.
\end{exa}

\subsection{Functional calculus}\label{sec-funct}
Let
$T=R_{\phi}\dvtx f\mapsto f*\phi$ be a convolution operator with
a symmetric probability density $\phi\in L^2(G)$.
Consider a function $\psi\dvtx [0,2]\ra[0,2]$ that is
increasing, continuous with continuous derivative and
which satisfies $\psi(0)=0$, $\psi(1)=1$, $\psi(2)<2$.
With such a function we associate
the operator
\[
\psi(I-T)= \int_0^2\psi(s)\,dE^{I-T}_s
\]
and
%
\begin{equation}\label{Tpsi}
T_\psi= I-\psi(I-T),\qquad T=R_\phi.
\end{equation}

\begin{lem} \label{lem-Tpsi1}
Let $\phi\in L^2(G)$ be a symmetric probability density. Let
$\psi\dvtx[0, 2]\ra[0,2]$ satisfies $\psi(0)=0$, $\psi(1)=1$,
$\psi(2)<2$ and assume that $\psi$ is increasing and continuous with
continuous derivative. Then $T_\psi$ defined at~(\ref{Tpsi}) is in
$V(G)$, and $T^{n}_\psi$ has finite trace for all $n\ge2$.
Further,\vspace*{-1pt} if $\phi=\xi*\xi$ with $\xi=\check{\xi}\in
L^2(G)\cap L^1(G)$, then $T=R_{\phi_\psi}$ with
$\phi_\psi=(\phi_\psi){\check{} } \in L^2(G)$ and $R_{\phi_\psi}$
bounded on $L^2(G)$.
\end{lem}
\begin{pf}
Note that the operators $\psi(I-T)$ and $T_\psi$ belong to the
von Neumann algebra $V(G)$. Further, from the elementary
fact that $|1-\psi(s)|\le C|1-s|$ on $[0,2]$, for some $C\in(0,\infty)$,
we deduce that\vadjust{\goodbreak}
$T^{2k}_\psi$ is a Hermitian nonnegative element in $V(G)$
which is dominated by
\[
CR^2_\phi=C^2\int_0^2 |1-s|^2\,dE^{I-T}_s.
\]
This last Hermitian element has finite trace equal to
$C\tau(R_\phi^2)=C\phi^{(2)}(e)$. Hence,
$T_\psi^{2k}$ has finite trace for $n\ge2$. This implies that
$T^{(2k+1)}_\psi$ has (extended) finite trace.

If $\phi=\xi*\xi$, then $T$ is Hermitian nonnegative and of finite trace.
Further $T_\psi$ is also Hermitian nonnegative and dominated by
$CT$. Hence $T_{\psi}=R_a^2$ with $a\in L^2(G)$, $R_a$
bounded on $L^2(G)$ and $\check{a}=a$. In particular,
$T_\psi=R_{\phi_\psi}$ with $\phi_\psi=a*a \in L^2(G)$.
This function is not a probability density, in general.
It is a probability density when $\psi$ is a Bernstein function; see,
for example,~\cite{BSCsubord}, Section 3.4, and~\cite{NJ}, Section 3.9.
\end{pf}
\begin{lem} \label{lem-Tpsi}
Let $\phi\in L^2(G)$ be a symmetric probability density such that
$\lim_{n\ra\infty}\phi^{(2n)}(e)= 0$.
Let $\psi\dvtx[0,2]\ra[0,2]$ be nonnegative
increasing, continuous with continuous derivative and
such that $\psi(0)=0$, $\psi(1)=1$, $\psi(2)<2$. Then the operator
$T_\psi\in V(G)$ defined at~(\ref{Tpsi}) is such that
$T^n_\psi$ has finite trace and, setting
%
\begin{equation}\label{Npsi}
N^\psi_\phi=N_\phi\circ\psi^{-1},
\end{equation}
we have
%
\begin{equation} \label{psiTn1}\quad
\tau(T^n_\psi)=\int_0^1 (1-s)^n\,dN^\psi_{\phi}(s) +O(a^n),\qquad a=\psi(2)-1
\in[0,1).
\end{equation}
\end{lem}
\begin{rem}
The hypothesis $\psi(2)<2$ insures that
the contribution coming from the spectrum of $I-R_\phi$
that lies in the interval
$(1,2)$ is exponentially small.
If $R_\phi$ is nonnegative [as a Hermitian operator on $L^2(G)$],
the value of $\psi$ in the interval $(1,2)$ becomes completely irrelevant
and
\[
\tau(T^n_\psi)=\int_0^1 (1-s)^n\,dN^\psi_{\phi}(s).
\]
\end{rem}
\begin{pf*}{Proof of Lemma~\ref{lem-Tpsi}}
Since $R^2_\phi=\int_0^2|1-\psi(s)|^2 \,dE^{I-T}_s$ has finite trace
equal to $\phi^{(2)}(e)$, the
the nondecreasing functions
$N_\phi(s)=\tau( E^{I-T}_s)$ [see definition~(\ref{def-Nphi})] and
$N^\sharp_\phi(s)=\tau(E^{I-T}_{(2-s,2)})$ are finite for all $s\in(0,1)$.
Further, since
$\phi^{(2n)}(e)\ra0$, we have $N_\phi(0)=0$
[i.e., there are no $L^2(G)$ solutions to $Tf=f$]. Hence,
\[
\tau(T_\psi^n)
= \int_0^1 \bigl(1-\psi(s)\bigr)^n \,d N_\phi(s)
+\int_0^1\bigl(1-\psi(2-s)\bigr)^n
\,dN^\sharp_\phi(s).
\]
The second integral is bounded by
\[
\biggl|\int_0^1\bigl(1-\psi(2-s)\bigr)^n \,dN^\sharp_\phi(s)\biggr| \le
\int_0^1|1-\psi(2-s)|^2\,dN^\sharp_\phi(s)
|\psi(2)-1|^{n-2}.
\]
Since $T^2=\int_0^2|1-\psi(s)|^2\,dE^{I-T}_s$ has finite trace $\phi^{(2)}(e)$
and $|1-\psi(2-s)|\le C|1-s|$, we obtain that
\[
\int_0^1|1-\psi(2-s)|^2\,dN^\sharp_\phi(s) \le C\phi^{(2)}(e).
\]

This yields the desired estimate since, by hypothesis,
$|\psi(2)-1|<1$.
\end{pf*}

To illustrate this lemma, we treat
the following simple test case.
\begin{theorem}\label{th-RRP}
Let $\phi\in L^2(G)$ be a symmetric positive probability density such that
\[
\phi^{(2n)}(e)\simeq n^{-D/2} \qquad\mbox{at infinity.}
\]
Let $\psi\dvtx[0,2]\ra[0,2]$ be nonnegative
increasing, continuous with continuous derivative and
such that $\psi(0)=0$, $\psi(1)=1$, $\psi(2)<2$. Assume further that
$\psi(s)\simeq(s/\ell(1/s))^\alpha$ at $0$, where $\alpha\in(0,\infty
)$ and
$\ell$ a positive function, slowly varying at infinity with de Bruijn conjugate
$\ell^{\#}$. Then
\[
\tau(T^n_{\psi})\simeq[n^{1/\alpha}\ell^{\#}(n^{1/\alpha})]^{-D/2}
\qquad\mbox{at infinity.}
\]
\end{theorem}
\begin{pf} This follows easily from Proposition~\ref{pro-NPhi2} and
Lemma~\ref{lem-Tpsi}, together with~\cite{BGT}, Proposition 1.5.15.
\end{pf}

Similar considerations, together with the arguments developed in
\cite{BPS}, Lem\-ma~2.3, Proposition 2.5, yield the following result
which is most useful when dealing with super-polynomial behaviors.
\begin{theorem} \label{th-RRR}
Let $\phi\in L^2(G)$ be a symmetric positive probability density.
Let $\pi\dvtx(0,\infty)\ra(0,\infty)$ be such that $\pi$ and
$t\mapsto t/\pi(t)$
are continuous increasing functions which tend to infinity at infinity.
Let $\psi\dvtx[0,2]\ra[0,2]$ be nonnegative
increasing, continuous with continuous derivative and
such that $\psi(0)=0$, \mbox{$\psi(1)=1$}, $\psi(2)<2$.
Set
%
\begin{equation}
\label{RRR}\pi_\psi^{-1}(t)=
t \psi^{-1}(1/t) \pi^{-1}\bigl(1/\psi^{-1}(1/t)\bigr).
\end{equation}

\begin{longlist}[(2)]
\item[(1)] Assume that there exists $c_1\in(0,\infty)$ such that,
for $n$ large enough,
\[
-\log\phi^{(2n)}(e)\ge c_1n/\pi(n).
\]
Then there exists $c_2\in(0,\infty)$ such that, for $n$ large enough,
\[
-\log\tau(T^{n}_\psi) \ge c_2 n/\pi_\psi(c_2 n) .
\]
\item[(2)] Assume that there exists $C_1\in(0,\infty)$ such that,
for $n$ large enough,
\[
-\log\phi^{(2n)}(e)\le C_1n/\pi(n).\vadjust{\goodbreak}
\]
Then there exists $C_2\in(0,\infty)$ such that, for $n$ large enough,
\[
-\log\tau(T^{n}_\psi) \le C_2 n/\pi_\psi(n/C_2) .
\]
\end{longlist}
\end{theorem}
\begin{pf} Let us observe that for a bijection $\pi$,
the two properties
(a) $\pi$ and $t\mapsto t/\pi(t)$ are increasing, and
(b) $t\mapsto\pi^{-1}(t)/t$
is increasing, are equivalent.
Further, given that $\psi$ is positive increasing,
property (a) for $\pi$ implies (b) for $\pi$ which implies (b) for $\pi
_\psi$
which finally implies (a) for $\pi_\psi$. The result now easily follows from
Proposition~\ref{pro-NPhi2} and Lemma~\ref{lem-Tpsi}.
\end{pf}

It is useful to illustrate Theorem~\ref{th-RRR} with some concrete examples.
Note that Theorem~\ref{th-RRR} allows us to treat upper and lower bounds
separately. For simplicity, we write down the examples in the context of
the rough equivalence $\simeq$.
\begin{exa} \label{exa-gamma-omega}
Assume that $-\log\phi^{(2n)}(e)\simeq\log n$ and that
$\psi(t)\simeq1/\break\ell(1/t)$ where $\ell$ is an increasing
slowly varying function tending to infinity at infinity and such that
\[
\log\ell^{-1}(t)\simeq t^\gamma\omega(t)^{1+\gamma},
\]
where $\gamma\in[0,\infty)$ and $\omega$ is a slowly varying function
at infinity with de Bruijn conjugate $\omega^{\#}$. Then
\[
-\log\tau(T^n_\psi)\simeq n^{\gamma/(1+\gamma)}
/\omega^{\#}\bigl(n^{1/(1+\gamma)}\bigr).
\]
\end{exa}
\begin{exa} \label{exa-gamma}Assume that
\[
-\log\phi^{(2n)}(e)\simeq n^\gamma,\qquad \gamma\in(0,1),
\]
and that
\[
\psi(t)\simeq t^\alpha/\ell(1/t),\qquad \alpha\in[0, \infty),
\]
where $\ell$ is an increasing
slowly varying function at infinity such that,
for every $a>0$, $\ell(t^a)\simeq\ell(t)$.
Then
\[
-\log\tau(T^n_\psi)\simeq [n/\ell(n)]^{\gamma_\alpha},\qquad \gamma_\alpha
=\frac{\gamma}{\gamma+\alpha(1-\gamma)}.
\]
\end{exa}
\begin{exa} \label{exa-gamma-l}
Assume that
\[
-\log\phi^{(2n)}(e) \le n/\pi(n)
\]
with $\pi$ positive increasing.
\begin{itemize}
\item Assume that $\pi(t)= t^{1-\gamma} \ell(t)$ with $\gamma\in(0,1]$
and $\ell$ slowly varying and satisfying $\ell(t^a)\simeq\ell(t)$
for all $a>0$. Then, for any $\alpha\in(0,1)$ and $\psi(t)= t^\alpha$,
we have
\[
-\log\tau(T^n_\psi)\le
[n/\ell(n)^{\alpha/\gamma}]^{\gamma_\alpha},\qquad
\gamma_\alpha
=\frac{\gamma}{\gamma+\alpha(1-\gamma)}.
\]
The cases $\gamma=1$ and $\gamma\in(0,1)$ should be treated
separately using slightly different arguments. See
\cite{BSCsubord}, Theorem 3.4,
for a similar computation.\vadjust{\goodbreak}
\item Assume that $\pi$ is regularly varying of index less than $1$.
Then for any positive increasing slowly varying $\ell$, $\psi=1/\ell(1/t)$,
and any $\varepsilon\in(0,1)$,
we have (see~\cite{BSCsubord}, Theorem 3.4,
for a similar computation)
\[
-\log\tau(T^n_\psi)\le C_\varepsilon n/\ell(\pi( C_\varepsilon n^\varepsilon)).
\]
\end{itemize}
\end{exa}

\subsection{Trace and comparison}
Let $T_1,T_2$ be
self-adjoint contractions that belong to a von Neumann algebra
$V$ equipped with a faithful semifinite normal trace~$\tau$.
For $i=1,2$, let $E^{I-T_i}_{s}$, $s\in[0,\infty)$,
be the (left-continuous) spectral projectors of $I-T_i$, so that
$T_i=\int_{0}^{\infty}(1-s)\,dE^{I-T_i}_s$. The following result is
crucial for our purpose. It is the von Neumann version
of a classical finite-dimensional spectral comparison theorem.
We set
\[
N_i(s)=\tau(E^{I-Ti}_s),\qquad s>0, i=1,2.
\]
Note that it can well be the case that $N_i(s)=\infty$.
\begin{pro} \label{pro-N21}
Referring to the above setting and notation,
let $T_1,T_2$ be
self-adjoint contractions that belong to the von Neumann algebra
$V$ equipped with a faithful semifinite normal trace $\tau$.
Assume that
\[
(I-T_1)\le C(I-T_2)
\]
and that $T_2$ is nonnegative.
Then we have
%
\begin{equation}
\forall s\in[0,1)\qquad N_2(s)\le N_1(Cs).
\end{equation}
\end{pro}
\begin{pf} Recall that, for any bounded self-adjoint operator
$S\in V$, $E^S_{(a,b)}$ denotes
the spectral projector associated to $S$ and the interval $(a,b)$.
By convention, the left-continuous
spectral resolution of $S$ is $E^S_s=E^S_{(-\infty,s)}$ so
that $S=\int_{-\infty}^{+\infty}s\,dE^S_s$ and
$E^S_{(a,b)}=\int_{(a,b)}\,dE^S_s$.
According\vspace*{1pt} to~\cite{BK}, Lemma 3, if $S_1,S_2$ are nonnegative
self-adjoint operators such that $S_2\le S_1$ then
(allowing for the possibility that the traces in question are infinite)
%
\begin{equation}\label{BK<}
\forall s>0\qquad \tau\bigl(E^{S_2}_{(s,\infty)}\bigr)
\le\tau\bigl(E^{S_1}_{(s,\infty)}\bigr).
\end{equation}
By hypothesis, we have $I-T_1\le C(I-T_2)$, which we write
\[
T_2\le I- C^{-1}(I-T_1).
\]
Applying~(\ref{BK<}) to $S_2=T_2$, $S_1= I-C^{-1}(I-T_1)$
($T_2$ is nonnegative by hypothesis and this implies that $S_2,S_1$
are also nonnegative)
and using the simple fact that
\[
E^{I-C^{-1}(I-T_1)}_{(s,\infty)}=E^{T_1}_{(1-C(1-s),\infty)},
\]
we obtain
\[
\forall s>0\qquad \tau\bigl(E^{T_2}_{(s,\infty)}\bigr)\le
\tau\bigl(E^{T_1}_{(1-C(1-s),\infty)}\bigr).
\]
Translating this inequality in terms of the spectral functions
\[
N_i(s)=\tau\bigl(E^{I-T_i}_{(-\infty,s)}\bigr)=
\tau\bigl(E^{T_i}_{(1-s,\infty)}\bigr),\vadjust{\goodbreak}
\]
we obtain
\[
\forall s\in[0,1)\qquad N_2(s)\le N_1(Cs).
\]
\upqed\end{pf}
\begin{cor}\label{cor-comp}
Referring to the above setting and notation, assume that $T_1,T_2$ are
nonnegative and that
there exist an integer $k_0$ and a constant $C\ge1$ such that
\[
\tau(T_1^{k_0}),\qquad\tau(T_2^{k_0})<\infty\quad\mbox{and}\quad
I -T_1\le C(I-T_2).
\]
Then, for all $n\ge k_0$,
\[
\tau(T_2^n)\le2C^2\tau\bigl(T_1^{\lfloor n/2C\rfloor}\bigr)+
2e^{-(n/16C)+k_0/8}\bigl(\tau(T^{k_0}_2)+
2C^2\tau(T^{k_0}_1)\bigr).
\]
\end{cor}
\begin{pf} We have
\begin{eqnarray*}
\tau(T_i^n) &=& \int_0^1(1-s)^n\,dN_i(s)\\
&=& n\int_0^\varepsilon(1-s)^{n-1}N_i(s)\,ds
+ (1-\varepsilon)^nN_i(\varepsilon)\\
&&{} + \int_\varepsilon^1(1-s)^n\,dN_i(s) .
\end{eqnarray*}
Since $(1-s)^{k_0}N_i(s)\le\tau(T^{k_0}_i)$
and $\int_0^1(1-s)^{k_0}\,dN_i(s)=\tau(T^{k_0}_i)$,
we obtain that
\[
\biggl|\tau(T^n_i)-n\int_0^\varepsilon(1-s)^{n-1}N_i(s)\,ds\biggr|\le
2(1-\varepsilon)^{n-k_0}\tau(T^{k_0}_i)
\]
for any real $n\ge k_0$.
Now, set $c=1/8C$, and use Proposition~\ref{pro-N21} and the elementary
inequality
$(1-s)\le(1-Cs)^{1/2C}$, $s\in[0,c]$, to write
\begin{eqnarray*}
n\int_0^{c}(1-s)^{n-1}N_2(s)\,ds&\le&
n\int_0^{c}(1-Cs)^{(n-1)/2C}N_1(Cs)\,ds\\
&\le& C n\int_0^{1/8}(1-s)^{(n-1)/2C}N_1(s)\,ds\\
&\le& 2C^2(n/2C) \int_0^{1/8}(1-s)^{(n/2C)-1}N_1(s)\,ds.\\
\end{eqnarray*}
It thus follows that
\[
\tau(T_2^n)\le2C^2 \tau(T_1^{n/2C}) +
2\biggl(1-\frac{1}{8C}\biggr)^{n-k_0}\tau(T^{k_0}_2)+
4C^2\biggl(1-\frac{1}{8}\biggr)^{(n/2C)-k_0}\tau(T^{k_0}_1).
\]
This yields the desired result.
\end{pf}

In applications of Corollary~\ref{cor-comp},
one may want to relax the hypothesis that $T_1,T_2$ are nonnegative.
This is possible thanks to the following result.\vadjust{\goodbreak}
\begin{cor}\label{cor-comp1}
Referring to the above setting and notation,
assume that there exist an integer $k_0$ and a constant $C\ge1$ such that
\[
\tau(T_1^{k_0}),\qquad\tau(T_2^{k_0})<\infty\quad\mbox{and}\quad
I -T_1\le C(I-T_2).
\]
Assume further that $\tau(T^k_2)\ge0$ for all $k\ge k_0$.
Then there are constants $C_1, C_2$ depending only on upper bounds on
$C, \tau(T^{k_0}_1), \tau(T_2^{k_0})$ and such that
\[
\tau(T_2^{2n})\le C_1\bigl(\tau\bigl(T_1^{2\lfloor n/C_2\rfloor}\bigr)+
e^{-n/C_2}\bigr) \qquad\mbox{for all $n$ large enough}.
\]
\end{cor}
\begin{pf}
Set $S=\frac{1}{2}(T_2^2+T_2^3)=\frac{1}{2}T_2^2(T_2+I)$.
This is a Hermitian nonnegative contraction. Further
$\tau(S^{2n})= 2^{-{2n}} \sum_0^{2n} {2n\choose i} \tau(T_2^{6n-i})$. Since
$\ell\mapsto\tau(T_2^{2\ell})$ is decreasing (e.g., by spectral theory)
and $\tau(T_2^{2\ell+1})\ge0$ (by hypothesis), we have
\[
\tau(S^{2n})\ge\frac{1}{2^{2n}}
\sum_{k\in2\mathbb N\cap[2n,6n]} \pmatrix{2n\cr k} \tau(T_2^{k})
\ge\frac{1}{2}\tau(T_2^{6n}).
\]
This shows that it suffices to estimate $\tau(S^{2n})$ by
$\tau(T_1^{2\lfloor c n\rfloor})$ for some $c>0$. This will follow
from Corollary~\ref{cor-comp} applied to the Hermitian nonnegative contractions
$T=T_1^2$, $S=\frac{1}{2}(T_2^3+T_2^2)$, if we can prove that
$I-T\le4C(I-S)$. This last inequality follows immediately from the hypothesis
$I-T_1\le C(I-T_2)$ because $I-T=I-T_1^2\le2( I-T_1)$ and
$I-T_2\le2(I-S)$. The last two inequalities follows from spectral theory
and the elementary inequalities $1-s^2\le2(1-s)$ and
$1-s\le2 -s^3-s^2$, $s\in[-1,1]$.
\end{pf}
\end{appendix}


%

%
\printaddresses

\end{document}